\documentclass[12pt]{amsart}
\usepackage[T1]{fontenc}
\usepackage[pdftex]{color,graphicx}
\usepackage{bbm} 
\usepackage{amsmath}
\usepackage{mathabx}
\usepackage[normalem]{ulem}
\usepackage{dsfont}
\usepackage{url}
\usepackage{mathtools}
\mathtoolsset{showonlyrefs}
\usepackage{anysize}
\usepackage{nicefrac}
\usepackage{geometry}
\marginsize{2.0cm}{2.0cm}{2.0cm}{2.0cm}
\usepackage{amssymb}
\usepackage{hyperref}
\usepackage{color}
\usepackage{orcidlink}
\theoremstyle{plain}
\numberwithin{equation}{section}
\newtheorem{theorem}{Theorem}[section]
\newtheorem{lemma}[theorem]{Lemma}
\newtheorem{corollary}[theorem]{Corollary}
\newtheorem{proposition}[theorem]{Proposition}
\newtheorem{remark}[theorem]{Remark}
\newtheorem{example}[theorem]{Example}
\newtheorem{hypothesis}[theorem]{Hypothesis}
\newtheorem{definition}[theorem]{Definition}
\newcommand{\ud}{\mathrm{d}}

\newcommand{\e}{\varepsilon}
\newcommand{\NN}{\mathbb{N}}
\newcommand{\RR}{\mathbb{R}}
\newcommand{\EE}{\mathbb{E}}
\newcommand{\gqq}{\geqslant}
\newcommand{\lqq}{\leqslant}
\newcommand{\ii}{\mathsf{i}}

\begin{document}
{
\begin{tikzpicture}[remember picture, overlay]
\node[anchor=north west, xshift=1.5cm, yshift=-1.8cm] at (current page.north west){
\includegraphics[width=4cm]{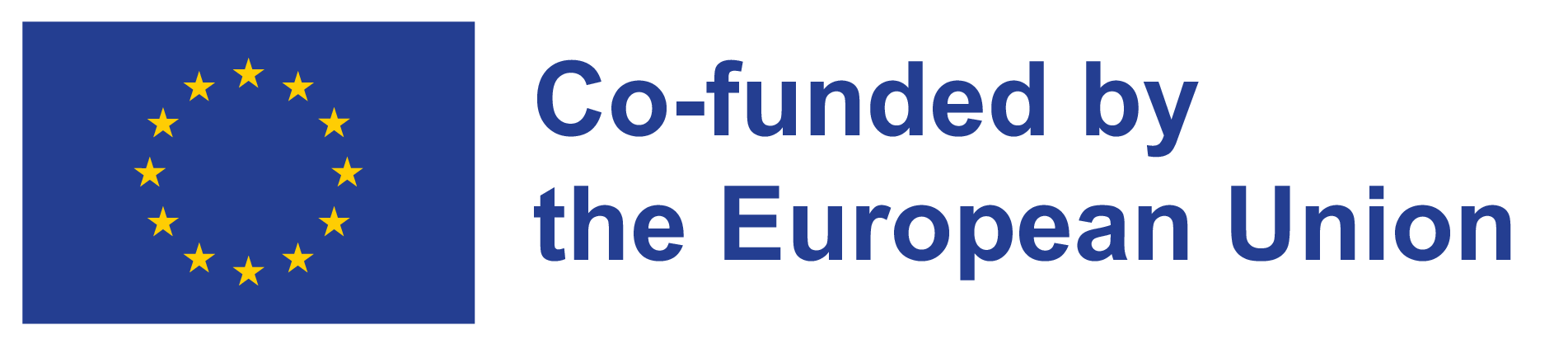} 
};
\end{tikzpicture}
}

\title[Cutoff stability for multivariate gBM with commutativities]{Cutoff stability of multivariate geometric Brownian motion}

\author{Gerardo~Barrera\orcidlink{0000-0002-8012-2600}}
\address{Center for Mathematical Analysis, Geometry and Dynamical Systems, Mathematics Department, Instituto
Superior T\'ecnico, Universidade de Lisboa, Lisboa,
Portugal.
\url{https://orcid.org/0000-0002-8012-2600}}
\email{gerardo.barrera.vargas@tecnico.ulisboa.pt}
\thanks{*Corresponding author: Gerardo Barrera.}

\author{Michael~A.~H\"ogele\orcidlink{0000-0001-5744-0494}}
\address{Departamento de Matem\'aticas, Facultad de Ciencias, Universidad de los Andes, Bogot\'a, Colombia.
\url{https://orcid.org/0000-0001-5744-0494}}
\email{ma.hoegele@uniandes.edu.co}

\author{Juan~Carlos~Pardo\orcidlink{0000-0003-4348-962X}}
\address{Center for Research in Mathematics CIMAT. Jalisco S/N, Valenciana, CP. 36240. Guanajuato, M\'exico.
\url{https://orcid.org/0000-0003-4348-962X}}
\email{jcpardo@cimat.mx}



\begin{abstract}
This article establishes cutoff convergence or abrupt convergence of three statistical quantities for multivariate (Hurwitz) stable geometric Brownian motion: the autocorrelation function, the Wasserstein distance between the current state and its degenerate limiting measure, and, finally, anti-concentration probabilities, which yield a fine-tuned trade-off between almost sure rates and the respective integrability of the random modulus of convergence using a quantitative Borel--Cantelli Lemma. We obtain in case of simultaneous diagonalizable drift and volatility matrices a complete representation of the mean square and derive nontrivial, sufficient and necessary mean square stability conditions, which include all real and imaginary parts of the volatility matrices' spectra.  
\end{abstract}
\maketitle
\section{\textbf{Introduction}}
Geometric Brownian motion and its multivariate generalizations is an important model class for many applications. It is applied in mathematical finance modeling of stock
prices and 
wealth dynamics
~\cite{Hu, LO20, MS10, MR05, PK13, PZ17,PZ16, PP06, Po21, Sh04,Stoj2012,Stojkoski},  in physics for turbulence theory~\cite{GCB23}, harmonic oscillators~\cite{Gi05, Gi05b}, 
inertia systems~\cite{HuHu18}, and hydro-dynamics~\cite{TE05,TTRD93}, and in biology and  medicine, modeling tumor 
stability~\cite{OOY16, SB05} and smoking habit dynamics~\cite{LOKB11, SG08}.
For more applications in diverse fields, we refer the introduction of~\cite{GCB23}.
It is also a well-studied mathematical object, 
since it is arguably one of the simplest systems 
where the difference in modeling with Stratonovich, It\^o or more general (H\"anggi-type) noises can be observed. 
Moreover, the linearity makes it suitable for calibration by model fitting, and for efficient simulation,
see~\cite{Appleby, Bycz, Dietz, Duf, Haskovec} and the references therein.
 
A geometric Brownian motion in one dimension is the solution of
the following linear stochastic differential equation (SDE) in the It\^o sense with linear multiplicative noise 
\begin{equation}\label{eq:BS}
\ud X_t(x)=aX_t(x)\ud t+bX_t(x)\ud W_t,\quad t\gqq 0,\quad X_0(x)=x,
\end{equation}
where $a$ and $b$ are given real numbers, $x\in \mathbb{R}$ is an initial datum, and $(W_t)_{t\gqq0}$ is a  scalar standard  Wiener process. The SDE~\eqref{eq:BS}
can be solved explicitly and has the following shape
\[
X_t(x)=\exp\left(\left(a-\frac{b^2}{2}\right)t+bW_t\right)x,\quad t\gqq 0,\quad x\in \mathbb{R},
\]
see for instance~\cite[Example~5.5]{Mao} and~\cite[Example~5.1.1]{Oks00}.
If~\eqref{eq:BS} is understood in the Stratonovich sense, where the correction term $\exp(-(\nicefrac{b^2}{2}) t)$ is contained in the stochastic (Stratonovich) integral, the solution has the shape $X_t(x)=\exp(at+bW_t)x$, $t\gqq 0$, $x\in \mathbb{R}$. This simple exponential form may change drastically in higher dimensions if $a$ and $b$ are replaced by deterministic square matrices $A\in \mathbb{R}^{d\times d}$ and $B\in \mathbb{R}^{d\times d}$, and it is only valid in special cases, see~\cite{Blanes,Hall} for deterministic setting.
In particular, a multivariate geometric Brownian motion in the Stratonovich sense
with initial value $x\in \RR^d$, $d\in \mathbb{N}$, has the following shape 
\begin{equation}\label{e:stochexp} 
X_t(x) = \exp(Y_t)x\quad \textrm{ with }\quad Y_t = At + BW_t,\qquad t\gqq0,   
\end{equation}
only if and only if the deterministic matrices $A$ and $B$ commute, see~\cite{Kamm,Yamato}.
By the Stratonovich change of variables formula, it is well-known that in case of commuting  matrices $A$ and $B$, the process $(X_t(x))_{t\gqq 0}$ given in~\eqref{e:stochexp} is the unique strong solution of the Stratonovich linear SDE 
\begin{equation}\label{eq:SDE}
\ud X_t(x) = A X_t(x)\ud t+ B  X_t(x)\circ \ud W_t, \quad t\gqq0, \quad  X_0(x)=x. 
\end{equation}
For non-commuting matrices $A$ and $B$, the solution of~\eqref{eq:SDE} has an exponential shape as~\eqref{e:stochexp}, however, due to 
the celebrated Baker--Campbell--Hausdorff--Dynkin formula, the
exponent process $(Y_t)_{t\geq 0}$ is more complicated and even nonlinear, and can be expressed in terms of iterated (nested) commutators of the Lie Algebra generated by $A$ and $B$, 
see~\cite{Kamm,Magnus,Muniz,Wang,Yamato,Yang}. We also remark that the exponential representation of the solution to~\eqref{eq:SDE} is challenging even in dimension two, see the example given by the so-called Ping-pong Lemma in~p.~1254 in~\cite{BBH}.
Due to its wide-spread use in applications, it is obvious that the long-term dynamics of (multivariate) geometric Brownian motion is of paramount importance.
Note that the dynamical equilibrium for the dynamics given in~\eqref{eq:SDE} is the null solution.
When the matrix $-A$ is (Hurwitz) stable, it is natural to expect that for any initial datum $x$, the 
$L^2$-convergence to the dynamical equilibrium  is exponentially fast. In the sequel, we connect the preceding convergence to the concept of cutoff phenomenon of Markov processes with small noise.

The concept of cutoff phenomenon was coined in the seminal papers~\cite{Al83,AD,AD87,DI,DS0,DS} by D.~Aldous, P.~Diaconis and
M.~Shahshahani
in the context of card-shuffling, random walks on the hypercube, and random walks on finite groups.
Since its discovery, the cutoff phenomenon has been rigorously established in numerous classes of Markov chains, spanning a broad range of probabilistic models and theoretical frameworks, see for instance~\cite{BY1,BASU,BD92,BBF08, CSC08,LanciaNardiScoppola12,LPW,MurrayPego} and the reference therein.
Recently, the concept of cutoff has been studied in machine learning for neural networks~\cite{AvelinKarlsson} and 
random quantum circuits~\cite{OhKais}. 
It studies a parametrized family of processes, a respectively parametrized family of limiting invariant measures and a respectively parametrized family of renormalized distances. In this setting, the cutoff phenomenon describes the situation of a critical time scale in terms of the  parameter (sometimes called complexity parameter) which separates sharply small and large values of the renormalized distance between the parametrized current state and the respective limiting measure. 

\textbf{The profile cutoff phenomenon:} 
For each $\e>0$, let $X^{\e}(x):=(X^\e_t(x_\e))_{t\gqq 0}$ be a stochastic process 
with values in the Polish space $E_\e$ (with at least two elements) and with initial position $x_\e\in E_\e$. 
Let $\mathcal{M}_1(E_\e)$ be the space of probability measures on $E_\e$ and consider $\widetilde{\mathcal{M}}_1(E_\e)$
 a subset of $\mathcal{M}_1(E_\e)$ 
equipped with the distance (or any way of measuring convergence) $\mathrm{dist}_\e$. 
Assume that the law of $X^\e_t(x_\e)$ belongs to $ \widetilde{\mathcal{M}}_1(E_\e)$ for each $\e>0$ and $t\geq 0$.
In addition, for each $\e>0$
there exists $\mu^\e \in \widetilde{\mathcal{M}}_1(E_\e)$ satisfying 
\begin{equation*}
\lim_{t\to \infty} \mathrm{dist}_\e(\mathrm{Law}(X^{\e}_t(x_\e)),\mu^\e)=0.
\end{equation*}
We say that a system $(X^{\e}(x_\e),\mu^\e,\mathrm{dist}_\e)_{\e>0}$
exhibits a \textit{profile cutoff phenomenon} at a cutoff time $t_\e\to \infty$ with cutoff window $w_\e=\mathrm{o}(t_\e)$, as $\e\to 0$,
if and only if the following limit exists 
\begin{equation}\label{eq:perfil}
\lim\limits_{\e\to 0}\mathrm{dist}_\e\big(\mathrm{Law}(X^\e_{t_\e+\rho\cdot w_\e}(x_\e)),\mu^\e\big)=\mathcal{P}(\rho) \quad \textrm{ for all } \quad \rho\in \mathbb{R},
\end{equation}
and additionally $\mathcal{P}(\infty)=0$ and $\mathcal{P}(-\infty)=D$,
where 
\[D:=\limsup\limits_{\e\to 0}\mathrm{Diameter}(\widetilde{\mathcal{M}}_1(E_\e),\mathrm{dist}_\e)\in (0,\infty].
\]
This means that the time scale $t_\e$ is a temporal threshold in the sense that surfing ahead of it $t_\e + \rho\cdot w_\e$, for $\rho>0$, yields small values while lagging behind the threshold sees maximal values of the respective limit, which is a function exclusively in terms of the deviation $\rho$ in $w_\e$ units from $t_\e$. 
In general, the limit~\eqref{eq:perfil} does not exist, and it is replaced by respective upper and lower limits, 
and which then gives rise to the weaker concept of \textit{window cutoff phenomenon}, see below. 
In the classical example of Aldous and Diaconis~\cite{AD, AD87} of card shuffling, $\e=1/n$ where $n$ represents the size of the deck of cards, $E_\e=S_n$ the space of card permutations (shufflings), $\mathrm{dist}_\epsilon$ is the usual non-normalized total variation distance, the diameter $D$ equals~$1$, and the profile $\mathcal{P}$ is given in terms of the Gaussian error function.
In general,
the question
whether a given family of uniquely ergodic stochastic processes exhibits a cutoff phenomenon is not straightforward. 
In addition, it may depend strongly of the underlying metrics between the probability measures, 
see~\cite{RachevKlebanov} for a structural discussion. 

\textbf{Profile cutoff convergence and its relation to asymptotic mixing:} Note that a profile cutoff phenomenon can occur even for systems, where the processes are not parameter dependent: $E_\e = E$, $x_\e = x$,  $X^\e = X$ and $\mu^\e = \mu$. 
In addition, we restrict ourselves to the cases where $\textrm{dist}_\e$ has the shape $\e^{-1} \cdot \textrm{dist}$ for some given parameter independent distance $\textrm{dist}$. For an example see the case study~\cite{BHOU} on the Ornstein--Uhlenbeck process with fixed noise amplitude and the Wasserstein distance. 
In the case of $\mathrm{dist}_\e=\mathrm{dist}/\e$ (diameter $D=\infty$) with a fixed distance $\mathrm{dist}$ on $\widetilde{\mathcal{M}}_1(E)$ the parameter $\e$ plays the role of an external parameter which quantifies the abrupt convergence of the non-normalized distance in the following sense: 
for small $\e>0$
\begin{equation}\label{eq:perfilstab}
\lim_{\e\to 0} \frac{\mathrm{dist}(\mathrm{Law}(X_{t_\e+\rho\cdot w_\e}(x)),\mu)}{\e} = \mathcal{P}(\rho)\quad 
\textrm{ for all } \quad \rho\in \mathbb{R}.
\end{equation}
We call this special type of degenerate profile cutoff phenomenon for non-parametrized systems \textit{profile cutoff stability}. This concept is linked naturally to an asymptotic version of $\gamma_\e$-mixing for small values $\e$. 
This connection can be outlined as follows. Assume that the map $\rho \mapsto \mathcal{P}(\rho)$ is monotonically decreasing and continuous, which can be guaranteed in full generality under Markovianity assumptions on the dynamics. 
Then, given the cutoff time scale $(t_\e)_{\e>0}$ and a cutoff window $(w_\e)_{\e>0}$, satisfying $w_\e \to w>0$, as $\e\to 0$ and a rate $\gamma_\e \searrow 0$, as $\e \to 0$, 
we can assume the existence of a deviation scale $\rho_\e$ satisfying
\begin{equation}\label{e:gammarhoscale}
\e \cdot \mathcal{P}(\rho_\e) = \gamma_\e \quad \textrm{ for all }\quad \e\in (0,1]. 
\end{equation}
In fact, if $\mbox{dist}$ denotes the Wasserstein distance of some order $p\gqq 1$, it turns out that for a variety of systems driven by additive L\'evy noise (see \cite{BHOU,BHPWA, BHPNO, BHPSPDE, BHPPESM}) the distance $\mbox{dist}(\mathrm{Law}(X_t(x)), \mu)$ has essentially the following shape 
$\mbox{dist}(\mathrm{Law}(X_t(x)), \mu) = \|x\| e^{-\lambda t} + \mathrm{o}(e^{-\lambda t})_{t\to \infty}$, 
where $\lambda$ is a positive spectral quantity. 
Then for $t_\e = \frac{1}{\lambda} \ln(1/\e)$, $w_\e \to w>0$, $\rho>0$ 
have in lowest order approximation (for short we write  $\dot =$) for the leading term 
\begin{equation}
\begin{split}
\frac{\mbox{dist}(\mathrm{Law}(X_{t_\e + \rho\cdot w_\e}(x)), \mu)}{\e} \dot = 
\|x\| \frac{e^{-\lambda (t_\e + \rho\cdot w_\e)}}{\e}
 = \|x\| e^{-\lambda \rho \cdot w_\e} \stackrel{\e\to 0}{\longrightarrow} \|x\| e^{-\lambda \rho \cdot w}=: \mathcal{P}(\rho).
\end{split}
\end{equation}
Therefore~\eqref{e:gammarhoscale} reads in this case as $\e e^{-\lambda \rho_\e \cdot w} = \gamma_\e$ and 
\begin{equation}
\frac{\mbox{dist}(\mathrm{Law}(X_{t_\e + \rho_\e \cdot w_\e}(x)), \mu)}{\gamma_\e} 
\dot =  \|x\| e^{-\lambda \rho_\e \cdot w_\e} \cdot \frac{\e}{\gamma_\e}  
= \|x\| e^{-\lambda \rho_\e (w_\e-w)}, 
\end{equation}
which yields the following: when the cutoff windows opens sufficiently fast in the sense that $\lim_{\e \to 0} \rho_\e (w_\e-w) = 0$ we obtain asymptotic $\gamma_\e$-mixing as follows:
\begin{equation}\label{asymptotic-mixing}
\lim_{\e \to 0} \frac{\mathrm{dist}(\mathrm{Law}(X_{t_\e+\rho_\e\cdot w_\e}(x)),\mu)}{ \|x\| \gamma_\e} = 1.
\end{equation}

\noindent \textbf{The window cutoff phenomenon and window cutoff convergence:} 
In the case that \[(\mathrm{dist}_\e\big(\mathrm{Law}(X^\e_{t_\e+\rho\cdot w_\e}(x_\e)),\mu^\e\big))_{\e>0}
\]
has more than one accumulation point, which may be generic for some systems (see for instance~\cite[Theorem~3.2]{BHPWA}), the natural generalization of the concept of profile cutoff phenomenon is the notion of the so-called \textit{window cutoff phenomenon}, that is,
\begin{equation}\label{eq:windowcutoff}
\begin{split}
\lim_{\rho\to \infty}\limsup\limits_{\e\to 0}\mathrm{dist}_\e\big(\mathrm{Law}(X^\e_{t_\e+\rho\cdot w_\e}(x_\e)),\mu^\e\big)&=0,\\
\lim_{\rho\to -\infty}\liminf\limits_{\e\to 0}\mathrm{dist}_\e\big(\mathrm{Law}(X^\e_{t_\e+\rho\cdot w_\e}(x_\e)),\mu^\e\big)&=D.\\
\end{split}
\end{equation}
In other words, the time scale $t_\e$ still splits large values from small values in the sense that a growing backward deviation measured in $w_\e$-units from $t_\e$ yields a maximal distance while 
a growing forward deviation measured in $w_\e$-units from $t_\e$ gives small values. 
In the analogous setting of profile cutoff stability~\eqref{eq:perfilstab} the notion of \textit{window cutoff stability} reads as follows:
there exist functions $\widehat{\mathcal{P}},\widecheck{\mathcal{P}}: \mathbb{R}\to [0,\infty)$ such that 
for small $\e>0$ 
\begin{equation}\label{d:pre-window-convergence}
 \begin{split}
&\lim_{\rho\to \infty}\limsup\limits_{\e\to 0} 
\frac{\mathrm{dist}(\mathrm{Law}(X_{t_\e+\rho\cdot w_\e}(x)),\mu)}{\e} 
= 
\lim_{\rho\to \infty} \widehat{\mathcal{P}}(\rho)=0,\\
&\lim_{\rho\to -\infty}\liminf\limits_{\e\to 0} 
\frac{\mathrm{dist}(\mathrm{Law}(X_{t_\e+\rho\cdot w_\e}(x)),\mu)}{\e} 
=
\lim_{\rho\to -\infty} \widecheck{\mathcal{P}}(\rho)= \infty.
\end{split} 
\end{equation}
We refer to the introductions~\cite{BHPWA,BHPTV, BJ1, BP,BY1,BLY06} for further details. 
In case that in~\eqref{eq:windowcutoff} and~\eqref{d:pre-window-convergence}, 
only the first relations are valid, respectively, then we (informally) speak of the upper pre-cutoff phenomenon or upper pre-cutoff convergence. 

\textbf{The main results:} 
The novelty of this article consists in the generalization of the concept of cutoff convergence to statistically insightful aspects of the threshold type dynamics of multivariate geometric Brownian motion $(X_t(x))_{t\gqq 0}$ defined as a solution of~\eqref{eq:SDE} with Hurwitz stable matrix $A$ and more general for the SDE~\eqref{eq:mstra} below.  
We show cutoff convergence for:
\begin{enumerate}
\item the autocorrelation matrix 
\[
R_t(x) =(\EE[(X_t(x))_{k} (X_t(x))_{k'}])_{k,k'\in\{1,\ldots,d\}},\] which is a fundamental statistical tool in time series analysis,
where $(X_t(x))_k$ denotes the $k$-th component of $X_t(x)$ for each $k\in \{1,\ldots,d\}$.
In Theorem~\ref{thm:auto} window cutoff convergence and in Theorem~\ref{thm:autoprof} profile cutoff convergence are shown under the minimal mean-square stability hypothesis of the system and without any algebraic (commutativity) assumption between the drift and noise coefficients. In particular, Theorem~\ref{thm:auto} implies (non-quantitative) $L^2$-window cutoff stability for the current state $\mathbb{E}[\|X_t(x)\|^2]$,
see Corollary~\ref{cor:norm}.
Since the dynamical equilibrium is degenerate $\delta_{0_d}$, the preceding $L^2$-convergence corresponds to the square of the so-called Wasserstein distance of order two.
\item the Wasserstein distance of order two between the law of the current state $X_t(x)$ and the degenerate limiting 
state $\delta_{0_d}$, that is,
\[
(\mathbb{E}[\|X_t(x)\|^2])^{1/2}.
\]
In Theorem~\ref{thm:solodiagonalizable} (quantitative) window cutoff convergence is shown under
simultaneously diagonalizable drift and noise coefficients. In addition,
a dynamical characterization of profile cutoff stability for simultaneously diagonalizable coefficients is given in Corollary~\ref{thm:profdiagonalizable}.
In particular,
(quantitative) profile cutoff stability for simultaneously unitarily diagonalizable coefficients is given explicitly in Corollary~\ref{thm:diagonalizable}.
We point out that in the case of simultaneously diagonalizable coefficients the imaginary parts of the eigenvalues of the dispersive matrices play role in the stability of~\eqref{eq:SDE}, see
Remark~\ref{rem:RI}, while for simultaneously unitarily diagonalizable coefficients do not play role.
\end{enumerate} 
Finally, we show upper pre-cutoff for 
\begin{enumerate}
\item[(3)] the anti-concentration probabilities 
\[
\mathbb{P}\left(\sup_{t\in (s_n, s_{n+1}]}\|X_t(x)\| >\e_n\right), \qquad \e_n>0
\]
for some increasing time scale $(s_n)_{n\in \NN}$, $s_n \to \infty$ in Proposition~\ref{thm:BCsupremum} under simultaneously unitarily diagonalizable drift and noise coefficients.
In combination with a quantitative version of the recently established Borel--Cantelli Lemma (Lemma~\ref{lem:BC}), this result in Corollary~\ref{cor:BCsupremum} yields a quantitative trade-off between the almost sure upper error bound $(\e_n)_{n\in \NN}$ 
and $(s_n)_{n\in \NN}$ such that 
\[
\sup_{t\in (s_n, s_{n+1}]}\|X_t(x)\| \lqq \e_n\quad \textrm{ for } \quad n > \mathcal{M}
\]
almost surely (a.s.) for the (random) modulus of convergence $\mathcal{M}$ and order of integrability of $\mathcal{M}$. 
For the special case of $s_n = t_n + \rho_n$, $\rho_n = \rho\cdot w_n$, $\rho, w_n>0$, $n\in \NN$, and $w_n \nearrow \infty$, as $n\to \infty$ we obtain an almost sure upper pre-cutoff trade-off between $(\rho_n)_{n\in \NN}$ and $(\e_n)_{n\in \NN}$. 
\end{enumerate}

\textbf{Additive vs. multiplicative noise:}
For linear and nonlinear SDEs perturbed by additive L\'evy noise, it is well-established in the literature~\cite{BHOU, BHPWA, BHPNO, BHPSPDE, BHPPESM}, that by the cost invariance property  of the Wasserstein distance of order $p\geq 1$ for the parallel transport (also referred to as as shift-linearity, see~\cite[Lemma~2.2(d)]{BHPWA}) the rate of the ergodic convergence, and more specifically the cutoff time, cutoff window, and cutoff profile, whenever it exists, do not depend of the noise characteristics. That is, only the spectral properties of the linear drift or the respective linearization of the nonlinear drift part at the stable point determine the speed of convergence to the limiting measure convergence. \\
\hfill
This decoupled behavior between the linear drift and the noise changes considerably in the presence of multiplicative noise. 
In~\cite{BHPSPDE}, Theorem~5.1 and Theorem~5.2, the authors establish a preliminary result with multiplicative noise and show the cutoff phenomenon for a stochastic partial differential equation with $\e$-small multiplicative noise. More precisely, they establish the cutoff phenomenon for the stochastic heat equation under Dirichlet conditions with $\e$-small amplitude multiplicative $Q$-Brownian motion and a Hilbert space valued L\'evy noise. In formula~(5.7) of~\cite{BHPSPDE} it can be seen that the multiplicative part enters the dynamics, however, only marginally, 
since it contains $\e$-small noise, such that the perturbation of the spectrum is also $\e$-small and hardly has any influence on the dynamics. \\
\hfill
In this manuscript, we establish the dependence of the ergodic convergence in finite dimensional linear equations with fixed (i.e. non-small) multiplicative noise on the noise coefficients. In particular, our results establish that under Stratonovich noise, $L^2$-convergence and cutoff-convergence can be established under Hurwitz stability of an associated enhanced matrix and its spectral properties. In order to study cutoff convergence in terms of the original system (not enhanced) we assume simultaneous diagonalizability of the drift and the "volatility" matrices. Moreover, we obtain an additional nontrivial, sufficient and necessary mean square stability condition (see Remark~\ref{rem:RI}), which is stronger, than only Hurwitz stability for the drift matrix, and includes additionally the real and imaginary parts of the spectra of the volatility matrices. In the last part of the article we show, how these results translate into almost sure error bounds which are uniform along a specific sequence of time intervals.

In the sequel, for clarity of the presentation, we recall some notation that it is used along the statements and proofs below.

\textbf{Notation:} The set of natural numbers is denoted by $\mathbb{N}:=\{1,2,\ldots,\}$.
For $\mathbb{K}\in \{\mathbb{R},\mathbb{C}\}$ and $d\in \mathbb{N}$ we denote by $\mathbb{K}^{d\times d}$ the set of $d\times d$ matrices with entries in $\mathbb{K}$, and $\textsf{GL}(\mathbb{K},d)$ denotes the set of $d\times d$ invertible matrices with entries in $\mathbb{K}$. 
$\mathrm{Re}(z)$, $\mathrm{Im}(z)$, $\overline{z}$ and $|z|$ denote respectively the real part, imaginary part, a complex conjugate and the complex modulus of
a given $z\in \mathbb{C}$.
Also, $\ii$ denotes the imaginary unit.
A matrix $U\in \RR^{d\times d}$ is called \textit{Hurwitz stable} ($U <0$, for short), if its spectrum (set of eigenvalues) $\mathsf{spec}(U) \subset \mathbb{C}_-:=\{z\in \mathbb{C}: \mathrm{Re}(z)<0\}$. 
 For $U\in \mathbb{K}^{d\times d}$ let $U^*$ be the \textit{conjugate transpose} (a.k.a. Hermitian transpose) of $U$, that is,
$(U^*)_{j,k}=\overline{U}_{k,j}$ for all $j,k\in \{1,\ldots,d\}$. 
In other words, $U^*=\overline{U}^T$, where $T$ denotes the classical transpose operator for matrices.
When $U=U^*$ we say that $U$ is a Hermitian matrix.
The set $U(d)$ denotes the  group of unitary matrices, that is, $U\in U(d)$ if and only if $U\in \textsf{GL}(\mathbb{C},d)$ and $U^*=U^{-1}$. 
Let $\mathbb{S}^1:=\{z\in \mathbb{C}:|z|=1\}$.
For $\mathbb{K}\in \{\mathbb{R},\mathbb{C},\mathbb{S}^1\}$ and $d\in \mathbb{N}$ we denote by $\textsf{Her}(\mathbb{K},d)$ the set of Hermitian matrices with coefficients in $\mathbb{K}$.
 For $U\in \mathbb{K}^{d\times d}$ let $\mathrm{Trace}(U)$ be the \textit{trace} of $U$, i.e.,
$\mathrm{Trace}(U):=\sum_{j=1}^{d}U_{j,j}$.
 For $\lambda_1,\ldots,\lambda_d\in \mathbb{K}$ let $\mathrm{diag}(\lambda_1,\ldots,\lambda_d)\in \mathbb{K}^{d\times d}$ be the diagonal matrix with 
\[(\mathrm{diag}(\lambda_1,\ldots,\lambda_d))_{j,j}=\lambda_j\quad \mathrm{ for }\quad j\in \{1,\ldots,d\}.\]
The \textit{Lie bracket} or \textit{commutator} is given by $[U, V] := UV -VU$ for $U, V\in \mathbb{K}^{d\times d}$.
 In a conscious abuse of notation, we denote  the zero element of $\mathbb{K}^{d\times d}$ and
$\mathbb{K}^{d^2\times d^2}$by $O$.  $I_d$ denotes the identity matrix in $\mathbb{K}^{d\times d}$.
$\otimes$ denotes the usual Kronecker product between matrices on $\mathbb{K}^{d\times d}$.

The manuscript is organized as follows. After the exposition of the common setup of Hurwitz stable multivariate geometric Brownian motion we show in three subsequent subsections three types of generalized cutoff convergence. In Subsection~\ref{ss:thprauto} we present a generalized window and profile cutoff convergence for the autoregression function with the help of Lemma~\ref{lem:jara} given in Appendix~\ref{ap:exp}. In Subsection~\ref{ss:thpr} we establish cutoff convergence (window and profile) for the $L^2$ distance under the additional assumption of simultaneous diagonalizability. In Subsection~\ref{ss:as} we describe an a.s. upper pre-cutoff trade-off relation under simultaneous unitary diagonalizability. In Section~\ref{sec:thprauto}, Section \ref{sec:thpr} and Section~\ref{a.s.} the results of Section~\ref{s:main} are shown in the respective order.

\section{\textbf{The setting and the main results}}\label{s:main}

\noindent 
In this section, we introduce the standard setting of this paper. 
We consider the unique strong solution $(X_t(x))_{t\gqq 0}$ of the following linear homogeneous stochastic differential equation (SDE for short) in the Stratonovich sense (here denoted by $\circ$)
on $\mathbb{R}^d$
\begin{equation}\label{eq:mstra}
\ud X_t(x) = AX_t(x)\ud t+\sum\limits_{j=1}^{L}B_jX_t(x)\circ\ud W_j(t),
\quad t\geq 0,\quad
X_0(x) = x\in \mathbb{R}^d,
\end{equation}
where $d,L\in \mathbb{N}$, $A,B_1,\ldots,B_L \in \mathbb{R}^{d\times d}$ and $W_j:=(W_j(t))_{t\geq 0}$ for $j\in \{1,\ldots,L\}$ are independent and identically distributed
standard one-dimensional Brownian motions. 
Since the SDE~\eqref{eq:mstra} has linear coefficients, the existence and uniqueness of the path-wise strong solution 
of~\eqref{eq:mstra} is established, for instance, by Theorem~3.1 in~\cite{Mao} or Theorem~4.5.3 in~\cite{Kloeden}. 
For stochastic integration in the Stratonovich sense (sometimes also known as Fisk--Stratonovich integral) 
we refer to \cite[Section II.7]{Pr04}. We denote by
$(\Omega, \mathcal{F},(\mathcal{F}_t)_{t\geq 0}, \mathbb{P})$ the underlying filtered complete probability space satisfying the usual conditions, see Definition~2.25 of~\cite{Karatzas}, where the independent standard Brownian motions $W_{j}$ for $j\in \{1,\ldots,L\}$ are defined, and denote by $\mathbb{E}$ the expectation with respect to the probability measure $\mathbb{P}$.

By the It\^o change of variable formula we have
\begin{equation}\label{eq:m}
\ud X_t(x)= \widetilde{A}X_t(x) \ud t+\sum\limits_{j=1}^{L}B_jX_t(x)\ud W_j(t),\quad t\geq 0,\quad
X_0(x)=x,
\end{equation}
where the drift matrix $\widetilde{A}$ is given by
\begin{equation}\label{ec:defAtilde}
\widetilde{A}:=A+\frac{1}{2}\sum\limits_{j=1}^{L}B^2_j.
\end{equation}
For details on the It\^o--Stratonovich change of variable we refer to Section~4.9 in~\cite{Kloeden}. 
In the sequel we establish three different types of cutoff convergence for geometric Brownian motion $(X_t)_{t\geq 0}$ given as the solution of \eqref{eq:m}. 

\bigskip 

\subsection{\textbf{Cutoff convergence for the autocorrelation function under Hurwitz stability}}\label{ss:thprauto}\hfill\\ 

\noindent In this section, we show the following. Under the assumption that the matrix
\begin{equation}
\label{eq:defL}
\mathbb{R}^{d^2\times d^2}\ni 
\Delta := \widetilde{A}\otimes I_d+ I_d\otimes \widetilde{A}+\sum_{j=1}^{L} B_j\otimes B_j 
\quad \textrm{is Hurwitz stable},
\end{equation}
where $\otimes$ denotes the usual Kronecker product between matrices on $\mathbb{R}^{d\times d}$,
we establish that the so-called autocorrelation matrix function associated to~\eqref{eq:mstra} has a window cutoff stability phenomenon, see Theorem~\ref{thm:auto} below. In addition, we prove that the  Euclidean norm of~\eqref{eq:mstra} exhibits a window cutoff stability phenomenon, see Corollary~\ref{cor:norm} below.

We recall that the autocorrelation matrix function $(R_t(x))_{t\geq 0}$, that is 
\begin{equation}\label{def:Rt}
R_t(x):=\mathbb{E}[X_t(x)(X_t(x))^*]=
\left(
\mathbb{E}[(X_t(x))_j(X_t(x))_k]\right)_{j,k\in \{1,\ldots,d\}},\quad \mathrm{ for }\quad t\geq 0,
\end{equation}
solves the following non-commutative matrix deterministic  differential equation
\begin{equation}\label{eq:auto}
\frac{\ud}{\ud t}R_t(x) =\widetilde{A}R_t(x)+R_t(x)\widetilde{A}^*+\sum\limits_{j=1}^{L}B_jR_t(x) B^*_j,\quad t\geq 0,
\quad
R_0(x)= xx^*,
\end{equation}
see for instance Theorem~3.2
in~\cite{Mao} or Theorem~(8.5.5) Item~b) in~\cite{Arnold}.
Using the linear isomorphism $\mathrm{vec}:\mathbb{R}^{d\times d}\to \mathbb{R}^{d^2}$ given by 
\[
\mathrm{vec}(V):=(V_{1,1},\ldots,V_{1,d},V_{2,1},\ldots,V_{2,d},\ldots,V_{d,1},\ldots,V_{d,d})
\]
for $V=(V_{i,j})_{i,j\in\{1,\ldots,d\}}\in \mathbb{R}^{d\times d}$,
we have the linear homogeneous matrix differential equation on $\mathbb{R}^{d^2}$
\begin{equation}\label{eq:autovect}
\frac{\ud }{\ud t} \mathrm{vec}(R_t(x)) =\Delta \mathrm{vec}(R_t(x)),\quad t\geq 0,\quad
\mathrm{vec}(R_0(x))= \mathrm{vec}(xx^*),
\end{equation}
where $\Delta$ is defined in~\eqref{eq:defL}.
In formula~(1.7) p.~1255 of~\cite{BBH} there is a minor typo in the definition of the operator $L$ there\footnote{The correct definition of $L$ in formula~(1.7) p.~1255 of~\cite{BBH} is $L=A\otimes I_n+ I_n\otimes A+\sum_{j=1}^{\ell} B_j\otimes B_j$}. For convenience we provide the proof of~\eqref{eq:autovect}. 
Applying the linear isormorphism $\mathrm{vec}$ to~\eqref{eq:auto} on both sides
and recalling that
$\mathrm{vec}(M_1RM_2)=(M^T_2\otimes M_1) \mathrm{vec}(R)$ for $M_1,M_2,R\in \mathbb{R}^{d\times d}$, see for instance, Lemma~\ref{lem:kron} in Appendix~\ref{ap:kron},
we obtain 
\begin{equation}\label{eq:autote}
\begin{split}
\frac{\ud}{\ud t}\mathrm{vec}(R_t(x)) &=\mathrm{vec}(\widetilde{A}R_t(x))+\mathrm{vec}(R_t(x)\widetilde{A}^*)+\sum\limits_{j=1}^{L}\mathrm{vec}(B_jR_t(x) B^*_j)\\
&=(I_d\otimes \widetilde{A})\mathrm{vec}(R_t(x))+
(\widetilde{A}\otimes I_d)
\mathrm{vec}(R_t(x))+\sum\limits_{j=1}^{L}
(B_j\otimes B_j)
\mathrm{vec}(R_t(x))\\
&=\Delta \mathrm{vec}(R_t(x))
\end{split}
\end{equation}
for all $t\geq 0$ with initial condition $\mathrm{vec}(R_0(x))= \mathrm{vec}(xx^*)$,
where we have used $\widetilde{A}^*=\widetilde{A}^T$ and $B^*_j=B^T_j$ for all $j\in \{1,\ldots,L\}$.
The solution of~\eqref{eq:autovect} is given by
\begin{equation}\label{eq:solexp1}
\mathrm{vec}(R_t(x))=
e^{\Delta t}\mathrm{vec}(xx^*)\quad \mathrm{ for }\quad t\geq 0,\, x\in \mathbb{R}^d.
\end{equation}
We point out 
that differential equations~\eqref{eq:auto}~and~\eqref{eq:autovect} are always valid and well-defined, regardless of any assumption on the spectrum for the matrix
 $\Delta$.
In what follows, we assume the following stability hypothesis for $\Delta$.
\begin{hypothesis}[Hurwitz stability]\label{hyp:Lest}
\hfill

\noindent
The matrix $\Delta$ defined in~\eqref{eq:defL} satisfies $\Delta<0$. 
\end{hypothesis}

\begin{remark}
We point out that the set of Hurwitz matrices of any fixed dimension $\mathcal{H}$ forms a (non-convex) cone, that is, for any $\alpha>0$ and $H\in \mathcal{H}$ we have $(\alpha H)\in \mathcal{H}$. 

It is worth noting  that 
Hurwitz stability of a matrix $M\in \mathbb{R}^{d\times d}$ is a strictly weaker condition than the coercivity of $-M$, that is, there exists a constant $\delta>0$ such that 
\[
x^TMx\leq  -\delta \|x\|^2_{\mathbb{R}^d}
\quad \textrm{ for all } \quad x\in \mathbb{R}^d.
\]
For instance, the matrix 
\[
M=\begin{pmatrix}
0 & -1\\
\lambda & -\lambda 
\end{pmatrix}
\]
for $\lambda \in (0,1/2)$
has two different eigenvalues $-\lambda/2\pm \ii \theta$ with $\theta:=\sqrt{\lambda(4-\lambda)}/2>0$, which implies that $M$ is Hurwitz and diagonalizable.
However,
it is not coercive.
Indeed, for $x=(x_1,x_2)^T$ we have 
$x^TMx=(1-\lambda)x_1x_2-\lambda x^2_2$ and taking $x_1=-x_2\neq 0$ we obtain
 $\delta\leq \lambda-1/2<0$, which contradicts coercivity.
\end{remark}

\begin{remark}[Mean-square asymptotically stability and asymptotically stability of the autocorrelation function]
\hfill

\noindent
We recall that the deterministic linear homogeneous matrix differential equation~\eqref{eq:autovect} is asymptotically stable
if and only if
for any initial datum $z\in \mathbb{R}^{d^2}$ it follows that
\[
\lim_{t\to \infty}\|e^{\Delta t}z\|_{\mathbb{R}^{d^2}}=0,
\]
which is equivalent to
$\Delta<0$,
see for instance~Theorem~8.2~in~\cite{Hespanhabook}.
We also recall that the SDE~\eqref{eq:mstra} is mean-square asymptotically exponentially stable if and only if 
for all $x\in \mathbb{R}^d$ there exist positive constants $C$ and $\gamma$ such that 
\begin{equation}\label{eq:Cgamma}
\mathbb{E}[\|X_t(x)\|^2_{\mathbb{R}^{d}}]\leq Ce^{-\gamma t}\|x\|^2_{\mathbb{R}^{d}}\quad \textrm{ for all }\quad t\geq 0.
\end{equation}
Due to the linearity of~\eqref{eq:mstra},
one can see that the mean-square exponentially asymptotically stability of~\eqref{eq:mstra} is equivalent to $\Delta<0$,
see Remark~(11.3.3) Item c) in~\cite{Arnold}.
For further details we refer to~\cite{Buckwar,GieslHafstein,HafsteinGudmundsson,Khasminskii}.
\end{remark}

In the sequel, using Lemma~\ref{lem:jara} in Appendix~\ref{ap:exp} we study the correct 
 rate of convergence for the long-term dynamics given in~\eqref{eq:solexp1}
when $\Delta<0$.
For $\Delta\in \RR^{d^2\times d^2}$, $\Delta < 0$ 
and a given $xx^*\neq O$,  Lemma~\ref{lem:jara} yields the existence of 
\begin{itemize}
\item[(1)] ${q}:={q}(xx^*)>0$,
\item[(2)] ${\ell}:={\ell}(xx^*), {m}:={m}(xx^*)  \in  \{1,\ldots, d^2\}$,
\item[(3)] ${\theta}_1:={\theta}_1(xx^*),\dots,{\theta}_{m}:={\theta}_m(xx^*) \in \RR$,
\item[(4)] and linearly independent vectors ${v}_1:={v}_1(xx^*),\dots,{v}_{m}:={v}_m(xx^*) \in  \mathbb{C}^{d^2}$
\end{itemize}
satisfying
\begin{equation}\label{eq:expdelta}
\lim_{t \to \infty} \left\|\frac{e^{{q} t}}{t^{{\ell}-1}} \exp(t\Delta)\mathrm{vec}(xx^*) - \sum_{k=1}^{{m}} e^{\ii  t{\theta}_k} {v}_k\right\|_{\mathbb{C}^{d^2}} =0.
\end{equation}
Moreover, there are positive constants ${K}_0 := {K}_0(xx^*)$ and  ${K}_1 :={K}_1(xx^*)$ such that 
\begin{equation}\label{eq:belowaboveti}
{K}_0\lqq \liminf_{t\rightarrow \infty}\left\|\sum_{k=1}^{{m}} e^{\ii  t{\theta}_k}{v}_k\right\|_{\mathbb{C}^{d^2}}
\lqq 
\limsup_{t\rightarrow \infty}\left\|\sum_{k=1}^{{m}} e^{\ii  t{\theta}_k}{v}_k\right\|_{\mathbb{C}^{d^2}}\lqq {K}_1.
\end{equation}
In particular,
\begin{equation}\label{ec:uno1}
\lim_{t \to \infty} \left|\left\|\frac{e^{{q} t}}{t^{{\ell}-1}} \exp(t\Delta)\mathrm{vec}(xx^*)\right\|_{\mathbb{C}^{d^2}}
-\left\|\sum_{k=1}^{{m}} e^{\ii  t{\theta}_k}{v}_k\right\|_{\mathbb{C}^{d^2}}\right| =0,
\end{equation}
\begin{equation}\label{eq:limiteinf}
\liminf_{t \to \infty} \left\|\frac{e^{{q} t}}{t^{{\ell}-1}} \exp(t\Delta)\mathrm{vec}(xx^*)\right\|_{\mathbb{C}^{d^2}}
=\liminf_{t \to \infty} 
\left\|\sum_{k=1}^{{m}} e^{\ii  t{\theta}_k}{v}_k\right\|_{\mathbb{C}^{d^2}},
\end{equation}
and
\begin{equation}\label{eq:limitesup}
\limsup_{t \to \infty} \left\|\frac{e^{{q} t}}{t^{{\ell}-1}} \exp(t\Delta)\mathrm{vec}(xx^*)\right\|_{\mathbb{C}^{d^2}}
=\limsup_{t \to \infty} 
\left\|\sum_{k=1}^{{m}} e^{\ii  t{\theta}_k} {v}_k\right\|_{\mathbb{C}^{d^2}}.
\end{equation}
Since the eigenvalues of $\Delta\in \mathbb{R}^{d^2\times d^2}$ come in pairs of complex conjugates, the vector
$\sum_{k=1}^{m} e^{\ii  t\theta_k} v_k$ belongs on $\mathbb{R}^{d^2}$ and we can indistinctly use $\|\cdot\|_{\mathbb{R}^{d^2}}$ and $\|\cdot\|_{\mathbb{C}^{d^2}}$. 
Using the preceding notation, we define the $\e$-mixing time as follows
\begin{equation}\label{eq:mixtime}
t_\e:=\frac{1}{{q}}\ln\left(\frac{1}{\e}\right)+\frac{{\ell}-1}{{q}}\ln\left(\ln\left(\frac{1}{\e}\right)\right),\qquad \e\in (0,1).
\end{equation}
We observe that the Hilbert--Schmidt norm (a.k.a. Frobenius norm) of $R_t(x)$ satisfies
\begin{equation}\label{eq:normrelation}
\begin{split}
\|R_t(x)\|_{\mathrm{HS}}:&=\left(\sum_{j,k=1}^d (R_t(x))^2_{j,k}\right)^{1/2}=\left(\|\mathrm{vec}(R_t(x))\|^2_{\mathbb{R}^{d^2}}\right)^{1/2}=\|e^{\Delta t}\mathrm{vec}(xx^*)\|_{\mathbb{R}^{d^2}}.
\end{split}
\end{equation}
We point out that no commutativity relations between the matrices $A,B_1,\ldots,B_L$ are assumed. 
The first result establishes abrupt convergence of the autocorrelation function under only the stability assumption on $\Delta$.  

\begin{theorem}[Window cutoff stability for the autocorrelation function]\label{thm:auto}\hfill

\noindent
Let $(X_t(x))_{t\gqq 0}$ be the unique strong solution of~\eqref{eq:mstra} and let $x\neq 0_d$.
Assume that the Hypothesis~\ref{hyp:Lest} is satisfied. Let $(t_\e)_{\e>0}$ be a function defined in~\eqref{eq:mixtime}
and $(w_\e)_{\e>0}$ be any function of positive numbers such that 
$w_\e\to w$ as $\e\to 0$, with $w>0$. For any $r\in \mathbb{R}$ it follows that
\begin{equation}\label{ec:upper}
\lim\limits_{r\to \infty}\limsup_{\e\to 0}\frac{1}{\e}\|R_{t_\e+r\cdot w_\e}(x)\|_{\mathrm{HS}}=0
\end{equation}
and
\begin{equation}\label{ec:lower}
\lim\limits_{r\to -\infty}\liminf_{\e\to 0}\frac{1}{\e}\|R_{t_\e+r\cdot w_\e}(x)\|_{\mathrm{HS}}=\infty,
\end{equation}
where $(R_t(x))_{t\gqq 0}$ is defined in~\eqref{def:Rt}.
\end{theorem}

The proof is found in Section~\ref{sec:thprauto}.

The second result provides a dynamical characterization of profile cutoff convergence, which generically link to some orthonormality properties of the vectors $\{{v}_1,\ldots,{v}_{{m}}\}$, see Remark~\ref{rem:orthorati} for further explanation.

\begin{theorem}[Profile cutoff stability for the autocorrelation function]\label{thm:autoprof}\hfill

\noindent
Let the hypotheses and notation of Theorem~\ref{thm:auto} be satisfied. For any $r\in \mathbb{R}$ it follows that
\begin{equation}\label{ec:lowerpro}
\lim_{\e\to 0}\frac{1}{\e}\|R_{t_\e+r\cdot w_\e}(x)\|_{\mathrm{HS}}={q}^{1-{\ell}}e^{-r{q}w}\|v(x)\|_{\mathbb{C}^{d^2}}
\end{equation}
for some $v(x)\in \omega(x)$
if and only if the map
\begin{equation}\label{eq:constant}
\omega(x)\ni x\mapsto \|v(x)\|_{\mathbb{C}^{d^2}}\quad \textrm{ is constant},
\end{equation}
where 
\begin{equation}
\omega(x):=\left\{v\in \mathbb{R}^{d^2}:
\textrm{ $v$ is an accumulation point of $\left(\sum_{k=1}^{{m}} e^{\ii  t{\theta}_k}{v}_k\right)_{t\geq 0}$}
\right\}.
\end{equation}
\end{theorem}

The proof is given in Section~\ref{sec:thprauto}.

\begin{remark}[Orthogonality and rationally independent]\label{rem:orthorati}\hfill

\noindent
We point out that an explicit characterization of $\omega(x)$ may not be easy to obtain, and hence the verification of~\eqref{eq:constant}can be difficult. 
In Theorem~3.2 of~\cite{BHPWA} it is shown an alternative characterization  of ~\eqref{eq:constant} in terms of some orthogonality when ``the angles'' satisfy the generic condition of being rationally independent.
\end{remark}

The following corollary implies the $L^2$-window cutoff stability for the current state $X_t(x)$ to its dynamical equilibrium $\delta_{0_d}$.

\begin{corollary}
[Window cutoff stability]\label{cor:norm}\hfill

\noindent
Let the hypotheses and notation of Theorem~\ref{thm:auto} be satisfied. For any $r\in \mathbb{R}$ it follows that
\begin{equation}
\begin{split}
\lim_{r \to \infty} \limsup_{\e\to 0} \frac{\EE\left[\|X_{t_\e + r \cdot w_\e}(x)\|^2_{\mathbb{R}^d}\right]}{\e} = 0
\end{split}
\end{equation}
and
\begin{equation}
\begin{split}
\lim_{r \to -\infty} \liminf_{\e\to 0} \frac{\EE\left[\|X_{t_\e + r \cdot w_\e}(x)\|^2_{\mathbb{R}^d}\right]}{\e} = \infty.
\end{split}
\end{equation}
\end{corollary}

The proof is found in Section~\ref{sec:thprauto}.

\begin{remark}[No closed form for the second moments]
\hfill

\noindent
For all $t\geq 0$ and $x\in \mathbb{R}^d$ we note that
\[
\|X_t(x)\|^2_{\mathbb{R}^d}=\mathrm{Trace}(R_t(x))=\sum_{j=1}^{d} (R_t(x))_{j,j}.
\]
Applying It\^o's formula 
in~\eqref{eq:mstra} for the smooth observable $\mathbb{R}^d\ni z\mapsto \|z\|^2_{\mathbb{R}^d}$, one can see that in general the SDE is not close, that is to say, the drift and the diffusion coefficients cannot be written as a function of $\|X_t(x)\|^2_{\mathbb{R}^d}$, see Remark~(8.5.8)~in~Chapter~8~of~\cite{Arnold}.
\end{remark}

\subsection{\textbf{Cutoff convergence for simultaneously diagonalizable coefficients}}\label{ss:thpr}\hfill\\ 

\noindent
In this section,
when the coefficients of~\eqref{eq:mstra} are
simultaneously diagonalizable and the system~\eqref{eq:mstra} is mean-square asymptotically stable, 
we establish the window cutoff convergence. In addition, we provide a dynamical characterization of profile cutoff convergence.
When the coefficients of~\eqref{eq:mstra} are simultaneously unitarily diagonalizable, we show profile cutoff stability and compute explicitly the profile function.

Recall that $L\in \NN$ and $A, B_1, \dots B_L\in \mathbb{R}^{d\times d}$ are the coefficient matrices of~\eqref{eq:SDE}. We start recalling the following basic definitions from linear algebra. 

\begin{definition}[Simultaneously diagonalizable and simultaneously unitarily diagonalizable]\label{def:sdud}
\hfill

\noindent
The set $\{A, B_1, \dots, B_L\}$ is called \textit{simultaneously diagonalizable}, 
if there is $P\in \textsf{GL}(\mathbb{C},d)$ such that 
\begin{equation}
\label{eq:simdia}
\begin{split}
D_A:&= P^{-1} A P=\mathrm{diag}(\alpha_1,\ldots,\alpha_d),\\
D_{B_j}:&= P^{-1} B_j P=\mathrm{diag}(\beta_{1,j},\ldots,\beta_{d,j})
\quad \textrm{ for all }\quad j\in \{1,\ldots,L\}
\end{split}
\end{equation}
are diagonal matrices with possible complex diagonal entries.

The set $\{A,B_1,\ldots,B_L\}$ is called \textit{simultaneously unitarily diagonalizable}, 
if there is $U$ in the unitary group $U(d)$,  such that 
\begin{equation}
\label{eq:simdiauni}
\begin{split}
D_A:&= U^{-1} A U=\mathrm{diag}(\alpha_1,\ldots,\alpha_d),\\
D_{B_j}:&= U^{-1} B_j U=\mathrm{diag}(\beta_{1,j},\ldots,\beta_{d,j})
\quad \textrm{ for all }\quad j\in \{1,\ldots,L\}
\end{split}
\end{equation}
are diagonal matrices with possible complex diagonal entries.
\end{definition}
For further details about simultaneous diagonalizability we refer to~\cite{Axler,Horn}.

Under simultaneous diagonizability of the coefficients 
we  obtain a full quantitative representation of the $L^2$ norm of $X$, 
see Remark~\ref{rem:autoexplanation} for more in the context.  

In the sequel,  we restrict our attention exclusively to norms on $\mathbb{R}^d$ and $\mathbb{C}^d$ and for simplicity on the exposition we omit the distinction between the norms.

\begin{lemma}
[Mean square representation for simultaneously diagonizable coefficients]
\label{lemma:diagonal}
\hfill

\noindent
Assume that the set of coefficients  $\{A,B_1,\ldots,B_L\}$ is simultaneously diagonalizable. 
Then for any $x\in \mathbb{R}^d$ and $t\geq 0$ it follows that
\begin{equation}
\label{eq:eqest}
\mathbb{E}[\|X_t(x)\|^2]=\sum_{i,k=1}^d \overline{y}_i\Gamma_{i,k}y_k \exp\left(t R_{i,k}\right)\exp\left(\ii t I_{i,k}\right),
\end{equation}
where $\Gamma:= (P^{-1})^* P^{-1}$, $y:=Px$,
and 
for $i,k\in \{1,\ldots,d\}$
\begin{equation}
\begin{split}
\mathbb{R}\ni  R_{i,k}:&=\mathrm{Re}\left(
\overline{\alpha}_i + \alpha_k  +\frac{1}{2}\sum_{j=1}^L (\overline{\beta}_{i,j} + \beta_{k,j})^2
\right),\\
\mathbb{R}\ni I_{i,k}:&=\mathrm{Im}\left(\overline{\alpha}_i + \alpha_k  +\frac{1}{2}\sum_{j=1}^L (\overline{\beta}_{i,j} + \beta_{k,j})^2
\right).
\end{split}
\end{equation}
\end{lemma}

The proof is given in Section~\ref{sec:thpr}.
\begin{remark}[Mean-square stability in terms of the spectra of the coefficient matrices]\label{rem:RI}
\hfill

\noindent
For $i,k\in \{1,\ldots,d\}$
we have
\begin{equation}
\begin{split}
R_{i,k}&=
\mathrm{Re}(\alpha_i)+\mathrm{Re}(\alpha_k)+\frac{1}{2}
\sum_{j=1}^{L} (\mathrm{Re}(\beta_{i,j})+\mathrm{Re}(\beta_{k,j}))^2-\frac{1}{2}\sum_{j=1}^{L} (-\mathrm{Im}(\beta_{i,j})+\mathrm{Im}(\beta_{k,j}))^2,\\
I_{i,k}&=
-\mathrm{Im}(\alpha_i)+\mathrm{Im}(\alpha_k)+
\sum_{j=1}^{L} (\mathrm{Re}(\beta_{i,j})+\mathrm{Re}(\beta_{k,j}))(-\mathrm{Im}(\beta_{i,j})+\mathrm{Im}(\beta_{k,j})).
\end{split}
\end{equation}
It is clear that $R_{i,k}=R_{k,i}$ and $I_{i,k}=-I_{k,i}$ for all $i,k\in \{1,\ldots,d\}$.
By~\eqref{eq:eqest} we obtain the nontrivial mean-square asymptotic stability condition
\begin{equation}\label{eq:meanst1}
\begin{split}
&\max_{1\lqq i\leq k \lqq d} R_{i,k}<0,
\end{split}
\end{equation}
which coincides with Condition~\eqref{eq:deltaneg} in Remark~\ref{rem:delta}.
\end{remark}

Taking advantage of Remark~\ref{rem:RI} we obtain the following mean-square stability result, which takes into account the precise dependence of the initial value. 

\begin{lemma}
[Mean square stability for simultaneously diagonalizable coefficients]
\label{lem:keep}
\hfill

\noindent
Keep the assumptions and the notation of Lemma~\ref{lemma:diagonal}. Additionally, we 
assume that 
\begin{equation}
\begin{split}
&\max_{1\lqq i\leq k \lqq d} R_{i,k}<0.
\end{split}
\end{equation}
Then for any $x\in \mathbb{R}^d$, $x\neq 0_d$, there is some $q:=q(x)>0$ such that
\[
0<
\liminf\limits_{t\to \infty}e^{q t}\mathbb{E}[\|X_t(x)\|^2]
\leq 
\limsup\limits_{t\to \infty}e^{q t}\mathbb{E}[\|X_t(x)\|^2]<\infty.
\]
In addition, we have that 
\begin{equation}\label{def:que}
q(x)=\min\{-R_{i,k}:(i,k)\in \{1,\ldots,d\}^2\quad \textrm{such that}\quad\overline{y}_i \Gamma_{i, k}y_k\neq 0\}.
\end{equation}
In particular, for 
\[
I(x):=\{(i,k)\in \{1,\ldots,d\}^2: -R_{i,k}=q\quad 
\textrm{such that}\quad \overline{y}_i \Gamma_{i, k}y_k\neq 0\}.
\]
it follows that
\begin{equation}\label{ec:limsupinff}
\begin{split}
\liminf\limits_{t\to \infty}e^{qt}\mathbb{E}[\|X_t(x)\|^2] 
&=2\sum_{(i,k)\in I;\,i\leq k}\mathrm{Re}\left(\overline{y}_i \Gamma_{i, k}y_k\exp(\mathsf{i} l_{i,k})\right)>0,\\
\limsup\limits_{t\to \infty}e^{qt}\mathbb{E}[\|X_t(x)\|^2] 
&=2\sum_{(i,k)\in I;\,i\leq k}\mathrm{Re}\left(\overline{y}_i \Gamma_{i, k}y_k\exp(\mathsf{i} u_{i,k})\right)<\infty
\end{split}
\end{equation}
for some $l_{i,k},u_{i,k}\in \mathbb{R}$.

\noindent Moreover, the limit $\lim\limits_{t\to \infty}e^{q t}\mathbb{E}[\|X_t(x)\|^2]$ exists if and only if 
\begin{equation}\label{eq:profcond}
\begin{split}
&\textsf{Her}(\mathbb{S}^1,\ell)\ni H\mapsto\sum_{(i,k)\in I}\overline{y}_i \Gamma_{i, k}y_k h_{i,k}\quad \textrm{ is constant},
\end{split}
\end{equation}
where $\ell:=\ell(x)=\# I$ and $H=(h_{i,k})_{i,k\in\{1,\ldots,\ell\}}\in \textsf{Her}(\mathbb{S}^1,\ell)$.
\end{lemma}

The proof is given in Section~\ref{sec:thpr}.

In case that the diagonalization can be carried out by a unitary change of coordinates we obtain that our results 
can be simplified significantly. In particular, the renormalized limit of the mean square always exists. 

\begin{corollary}[Mean-square representation for simultaneously \textit{unitarly}
diagonalizable coefficients]\label{cor:simuluni}\hfill

\noindent
Assume that the set of coefficients $\{A,B_1,\ldots,B_L\}$ is simultaneously unitarily diagonalizable. 
Then for any $x\in \mathbb{R}^d$ and $t\geq 0$ it follows that
\[
\mathbb{E}[\|X_t(x)\|^2]=\sum_{i=1}^d |y_i|^2 \exp\left(t R_{i}\right),
\]
where  $y:=Ux$,
and 
for $i\in \{1,\ldots,d\}$
\begin{equation}
R_i:=2\mathrm{Re}(\alpha_i)+2\sum_{j=1}^{L}(\mathrm{Re}(\beta_{i,j}))^2.
\end{equation}
In addition, assume that $x\neq 0_d$ and 
\begin{equation}\label{def:eq1}
q:=\min\{-R_i:i\in\{1,\ldots,d\}\textrm{ and } y_i\neq 0\}>0.
\end{equation}
For
$I:=\{i\in \{1,\ldots,d\}: -R_i=q \textrm{ and } y_i\neq 0 \}$ it follows that
\begin{equation}\label{eq:limitfp}
\lim\limits_{t\to\infty}e^{qt}\mathbb{E}[\|X_t(x)\|^2] =\sum_{i\in I} |y_i|^2>0.
\end{equation}
\end{corollary}

The result is a consequence of Remark~\ref{rem:RI} and Lemma~\ref{lem:keep}. 

The following result under simultaneous diagonizability is a sharpened and quantified version of 
Corollary~\ref{cor:norm}, which takes into account the dependence of the initial condition. 

\begin{theorem}[Window cutoff stability for simultaneously diagonalizable coefficients]\label{thm:solodiagonalizable}\hfill

\noindent
Keep the assumptions and notation of Lemma~\ref{lem:keep}.
For any $x\in \RR^d$, $x\neq 0_d$, we consider $q:=q(x)>0$ given in~\eqref{def:que}~of~Lemma~\ref{lem:keep}. 
Define the cutoff time $(t_\e)_{\e\in (0,1]}$ by
\begin{equation}
t_\e :=t_\e(x)=\frac{1}{q}\ln\left(\frac{1}{\e}\right), \quad \e\in (0,1], \quad 
\end{equation}
and the width of the time window by $(w_\e)_{\e\in (0,1]}$ such that $w_\e\to w>0$ as $\e\to 0$. Then the following  window cutoff convergence is valid
\begin{equation}\label{eq:thstaresult}
\begin{split}
\liminf_{\e\to 0} \frac{\EE[\|X_{t_\e + r \cdot w_\e}(x)\|^2]}{\e}&= 
e^{-qr\cdot w}\sum_{(i,k)\in I;\,i\leq k}\overline{y}_i \Gamma_{i, k}y_k\exp(\mathsf{i} l_{i,k})\\
&=e^{-qr\cdot w}\sum_{(i,k)\in I;\,i\leq k}2\,\mathrm{Re}\left(\overline{y}_i \Gamma_{i, k}y_k\exp(\mathsf{i} l_{i,k})\right)>0,\\
\limsup_{\e\to 0} \frac{\EE[\|X_{t_\e + r \cdot w_\e}(x)\|^2]}{\e}&= 
e^{-qr\cdot w}\sum_{(i,k)\in I;\,i\leq k}\overline{y}_i \Gamma_{i, k}y_k\exp(\mathsf{i} u_{i,k})\\
&=e^{-qr\cdot w}\sum_{(i,k)\in I;\,i\leq k}\,2\mathrm{Re}\left(\overline{y}_i \Gamma_{i, k}y_k\exp(\mathsf{i} u_{i,k})\right)<\infty
\end{split}
\end{equation}
for some $l_{i,k},u_{i,k}\in \mathbb{R}$.
\end{theorem}

The proof is given in Section~\ref{sec:thpr}.

\begin{remark}[Eigenvalues of $\Delta$ when the set $\{A,B_1,\ldots,B_L\}$ is simultaneously diagonalizable]\label{rem:delta}\hfill

\noindent
Recall the autocorrelation matrix $\Delta $ given in~\eqref{eq:defL}.
Let
$P\in \textsf{GL}(\mathbb{C},d)$ such that 
\begin{equation}
\begin{split}
D_A:&= P^{-1} A P=\mathrm{diag}(\alpha_1,\ldots,\alpha_d),\\
D_{B_j}:&= P^{-1} B_j P=\mathrm{diag}(\beta_{1,j},\ldots,\beta_{d,j})
\quad \textrm{ for all }\quad j\in \{1,\ldots,L\}.
\end{split}
\end{equation}
By~\eqref{ec:defAtilde} we have
\begin{equation}
D_{\widetilde{A}}:=P^{-1} \widetilde{A} P=\mathrm{diag}(\widetilde{\alpha}_1,\ldots,\widetilde{\alpha}_d),
\end{equation}
where
\begin{equation}
\widetilde{\alpha}_k:=\alpha_k+\frac{1}{2}\sum_{j=1}^{L}\beta^2_{k,j}, \quad k=1,\ldots,d.
\end{equation}
Recall the properties given in Lemma~\ref{lem:kron} in Appendix~\ref{ap:kron}.
Define $Q=P\otimes P$. Since $\textrm{det}(Q)=(\textrm{det}(P))^{2d}> 0$, we have $Q\in \textsf{GL}(\mathbb{C},d^2)$. 
We claim that $Q^{-1}\Delta Q$ is a diagonal matrix.
Indeed,
\begin{equation}
\begin{split}
Q^{-1}\Delta Q&=(P\otimes P)^{-1}\Delta (P\otimes P)=(P^{-1}\otimes P^{-1})\Delta (P\otimes P)\\
&=(P^{-1}\otimes P^{-1})\big(\widetilde{A}\otimes I_d+ I_d\otimes \widetilde{A}+\sum_{j=1}^{L} B_j\otimes B_j \big) (P\otimes P)\\
&=(P^{-1}\otimes P^{-1})(\widetilde{A}\otimes I_d)(P\otimes P)
+(P^{-1}\otimes P^{-1})(I_d\otimes \widetilde{A})(P\otimes P)\\
&\qquad+\sum_{j=1}^{L} (P^{-1}\otimes P^{-1})(B_j\otimes B_j ) (P\otimes P)\\
&=(P^{-1}\widetilde{A}P)\otimes I_d
+I_d \otimes (P^{-1}\widetilde{A}P)+\sum_{j=1}^{L} (P^{-1}B_jP)\otimes (P^{-1}B_jP)\\
&=D_{\widetilde{A}}\otimes I_d
+I_d \otimes D_{\widetilde{A}}+\sum_{j=1}^{L} D_{B_j}\otimes D_{B_j},
\end{split}
\end{equation}
which is a diagonal matrix. Hence, $Q$ diagonalizes $\Delta$ and its eigenvalues are given by
\begin{equation}
\lambda_{k,k'}:=\widetilde{\alpha}_k+\widetilde{\alpha}_{k'}+\sum_{j=1}^{L} \beta_{k,j}\beta_{k',j}
=\alpha_k+\alpha_{k'}+\frac{1}{2}\sum_{j=1}^{L} (\beta_{k,j}+\beta_{k',j})^2
\end{equation}
for $k,k'\in \{1,\ldots,d\}$.
As a consequence, Hypothesis~\ref{hyp:Lest}, that is, $\Delta<0$, is equivalent to
\begin{equation}\label{eq:deltaneg}
\max\limits_{k,k'\in \{1,\ldots,d\} } \mathsf{Re}(\lambda_{k,k'})<0.
\end{equation}
\end{remark}

\begin{remark}\label{rem:autoexplanation}~
\begin{enumerate}
\item For the case of simultaneously diagonalizable coefficient matrices $\{A,B_1,\ldots,B_L\}$, formula~\eqref{eq:thstaresult} implies the window cutoff convergence in Corollary~\ref{cor:norm}.
However, 
the result is much finer since it gives the full quantification of the limits as exponentials which resembles the cutoff profile in Corollary~\ref{thm:profdiagonalizable} except for  different constants for the upper and lower limit regime.
\item Corollary~\ref{cor:norm} only assumes Hurwitz stability for the Kronecker matrix $\Delta$ given in~\eqref{eq:defL} with cutoff time given in terms of spectral properties of $\Delta$, that is, 
only indirectly and nontrivially in terms of the spectral properties of $\{A,B_1,\ldots,B_L\}$.
In the special case of Theorem~\ref{thm:solodiagonalizable} those can be directly identified in terms of the spectral properties of $\{A,B_1,\ldots,B_L\}$, see Lemma~\ref{lem:keep}.
\item We point out that Corollary~\ref{cor:norm}  covers the full open (non-convex) cone of Hurwitz stable matrices $\Delta$, while Theorem~\ref{thm:solodiagonalizable} treats only a closed set of this cone that corresponds to simultaneously diagonalizable coefficients.
\end{enumerate}
\end{remark}

The following result gives an algebraic characterization, when the upper and the lower limit in Theorem~\ref{thm:solodiagonalizable} in fact coincide and form the cutoff profile. 

\begin{corollary}[Characterization of profile cutoff stability for simultaneously diagonizable coefficients]\label{thm:profdiagonalizable}\hfill

\noindent
Keep the assumptions and notation of Lemma~\ref{lem:keep}.
For any $x\in \RR^d$, $x\neq 0_d$, we consider $q:=q(x)>0$ given in~\eqref{def:que}~of~Lemma~\ref{lem:keep}. 
Define the cutoff time $(t_\e)_{\e\in (0,1]}$ by
\begin{equation}
t_\e :=t_\e(x)=\frac{1}{q}\ln\left(\frac{1}{\e}\right), \quad \e\in (0,1], \quad 
\end{equation}
and the width of the time window by $(w_\e)_{\e\in (0,1]}$ such that $w_\e\to w>0$ as $\e\to 0$. Then the following  profile cutoff convergence is valid
\begin{equation}
\begin{split}
\lim_{\e\to 0} \frac{\EE[\|X_{t_\e + r \cdot w_\e}(x)\|^2]}{\e} =e^{-qr\cdot w}\sum_{(i,k)\in I}\overline{y}_i \Gamma_{i, k}y_k h_{i,k}\quad \textrm{ for any } \quad r\in \mathbb{R},
\end{split}
\end{equation}
 where $\ell:=\ell(x)=\# I$ and any $H=(h_{i,k})_{i,k\in\{1,\ldots,\ell\}}\in \textsf{Her}(\mathbb{S}^1,\ell)$,
if and only if, 
the function
\[
\textsf{Her}(\mathbb{S}^1,\ell)\ni H\mapsto\sum_{(i,k)\in I}\overline{y}_i \Gamma_{i, k}y_k h_{i,k}\quad \textrm{ is constant}.
\]
\end{corollary}

The proof is given in Section~\ref{sec:thpr}.

In case of unitary diagonalization the cutoff profile can be simplified considerably as is stated below. 

\begin{corollary}[Profile cutoff stability for unitarily diagonalizable coefficients]\label{thm:diagonalizable}\hfill

\noindent
Keep the assumptions and notation of Corollary~\ref{cor:simuluni}.
For any $x\in \RR^d$, $x\neq 0_d$, we consider $q:=q(x)>0$ given in~\eqref{def:eq1}~of~Corollary~\ref{cor:simuluni}. 
Define the cutoff time $(t_\e)_{\e\in (0,1]}$ by
\begin{equation}
t_\e :=t_\e(x)=\frac{1}{q}\ln\left(\frac{1}{\e}\right), \quad \e\in (0,1], \quad 
\end{equation}
and the width of the time window by $(w_\e)_{\e\in (0,1]}$ such that $w_\e\to w>0$ as $\e\to 0$. Then the following  profile cutoff convergence is valid 
\begin{equation}
\begin{split}
\lim_{\e\to 0} \frac{\EE[\|X_{t_\e + r\cdot w_{\e}}(x)\|^2]}{\e} =e^{-qr w}\sum_{i\in I} |y_i|^2\quad \textrm{ for any }\quad r\in \mathbb{R}.
\end{split}
\end{equation}
\end{corollary}

The proof is given in Section~\ref{sec:thpr}.

\begin{remark}
This result can be compared to the profile cutoff result for the heat equation with multiplicative noise. 
In fact, is the finite-dimensional analogue of Theorem~5.1 for constant noise intensity in~\cite{BHPSPDE}.  
\end{remark}

\begin{remark}[Monotonic decay of the Wasserstein distance]\hfill

\noindent
We stress that the Dirac measure at zero, $\delta_{0_d}$, is invariant for the dynamics~\eqref{eq:SDE} and hence 
\[
\mathcal{W}^2_2(X_t(x),\delta_{0_d}):=\mathbb{E}[\|X_{t}(x)\|^2],\quad \textrm{for any}\quad t\geq 0, \quad x\in \mathbb{R}^d,
\]
where $\mathcal{W}_2$ is the standard Wasserstein distance of order $2$, see~\cite{villani}.
Moreover, the map 
\begin{equation}\label{eq:monotona}
t\mapsto \mathcal{W}_2(X_t(x),\delta_{0_d}) 
\end{equation}
is known to be non-increasing, see for instance Lemma~B.3 (Monotonicity) 
in~\cite{Boursier}.
\end{remark}
The following result connects the cutoff stability with the notion of mixing times with respect to the Wasserstein distance of order $2$ and the respective cutoff phenomenon in the sense of Levin, Peres and Wilmer given in Chapter~18 of~\cite{LPW}, see the definition in~(18.3).
For further details about Wasserstein distances we refer to~\cite{villani}. 

\begin{corollary}[Asymptotic $\e$-mixing time]\label{cor:mixingI}
\hfill

\noindent
Assume the hypotheses and notation of 
Theorem~\ref{thm:solodiagonalizable}.
Given $\delta>0$ we define the $\delta$-mixing time as follows.
\[
\tau^x_\e(\delta):=\inf\left\{t\gqq0: \frac{\mathbb{E}[\|X_{t}(x)\|^2]}{\e}\lqq \delta\right\}.
\] 
Then for any $M>\delta$ it follows that
\begin{equation}\label{eq:mezcla}
\lim\limits_{\e\to 0} \frac{\tau^x_\e(\delta)}{\tau^x_\e(M-\delta)}=1\quad \mathrm{ and }\quad 
\lim\limits_{\e\to 0}\frac{\tau^x_\e(\delta)}{t_\e}=1.
\end{equation}
\end{corollary}
The proof is given in Section~\ref{sec:thpr}.
\bigskip 

\subsection{\textbf{Almost sure uniform rates and almost sure upper pre-cutoff for simultaneously unitarily diagonalizable coefficients}}\label{ss:as}
\hfill

The aim of this section is to provide an a.s. quantification of cutoff convergence with the help of a recently established quantitative version of the first Borel--Cantelli Lemma
(see~\cite{EHS})
by making use of estimates of the second moment established 
in Section~\ref{a.s.}.

In this section, we assume that the set of coefficients
$\{A,B_1,\ldots,B_L\}$ is simultaneously unitarily diagonalizable.
We start noticing that $\{A,B_1,\ldots,B_L\}$ is simultaneously unitarily diagonalizable if and only if $\{A^*,B^*_1,\ldots,B^*_L\}$ is simultaneously unitarily diagonalizable.
Then there exists an ordered orthonormal basis $\{v_1, \dots, v_d\}$ of $\mathbb{C}^d$ such that 
\begin{equation}\label{eq:basis}
\begin{split}
A^* v_k &= \lambda_k v_k,\quad k\in \{1,2,\ldots,d\},\quad \textrm{ and }\quad\\
\qquad B_j^* v_k &= \mu_{k,j} v_k, \quad j\in \{1,2,\ldots, L\},\,k\in \{1,2,\ldots,d\}.
\end{split}
\end{equation}
Note that in general  $v_k\in \mathbb{C}^d$ for all $k\in \{1,2,\ldots,d\}$, and $\lambda_k,\mu_{k,j}\in \mathbb{C}$ for all $j\in \{1,2,\ldots, L\},\,k\in \{1,2,\ldots,d\}$. In fact, using the notation in~\eqref{eq:simdiauni} one can assume the labeling  $\lambda_k=\overline{\alpha}_k$
and $\mu_{k,j}=\overline{\beta}_{k,j}$ for all $j\in \{1,\ldots, L\}$ and $k\in \{1,\ldots,d\}$.
Furthermore, we set
\begin{equation}\label{def:tildelamb}
\begin{split}
\widetilde{\lambda}_{k} &= \lambda_k + \frac{1}{2} \sum_{j=1}^{L}\mu^2_{k,j}\quad  \textrm{ for all }\quad k\in \{1,2,\ldots,d\}.
\end{split}
\end{equation}
Note that for all $k\in \{1,\ldots,d\}$
\begin{align}
&\mathsf{Re}(\widetilde{\lambda}_{k})=
\mathsf{Re}\left(\lambda_{k}\right)+\frac{1}{2}\sum_{j=1}^{L}
\left((\mathsf{Re}(\mu_{k,j}))^2-(\mathsf{Im}(\mu_{k,j}))^2\right),\\
&\mathsf{Im}(\widetilde{\lambda}_{k})=
\mathsf{Im}\left(\lambda_{k}\right)+\sum_{j=1}^{L}
\mathsf{Re}(\mu_{k,j})\mathsf{Im}(\mu_{k,j})\quad \textrm{ and }\\
&\mathsf{Re}(\widetilde{\lambda}_k)  +  \frac{1}{2}\sum_{j=1}^L |\mu_{k,j}|^2=
\mathsf{Re}\left(\lambda_{k}\right)+\sum_{j=1}^{L}
(\mathsf{Re}(\mu_{k,j}))^2.\label{eq:ghy}
\end{align}
Then~\eqref{eq:basis} gives for all $k\in \{1,2,\ldots,d\}$ that
\begin{equation}\label{eq:veryeigen}
\widetilde{A}^* v_k =\left(A^* + \frac{1}{2} \sum_{j=1}^L (B_j^*)^2\right) v_k = \widetilde{\lambda}_k v_k. 
\end{equation}

In the following result, we take advantage on the fine structure introduced before to obtain tail estimates. 

\begin{proposition}[Quantitative tail estimate of the anti-concentration probabilities under mean-square stability]\label{thm:BCsupremum}
\hfill

\noindent
Assume that the set of coefficients
$\{A,B_1,\ldots,B_L\}$ is simultaneously unitarily diagonalizable and
\begin{equation}\label{e:uniformLambda}
\Lambda:=\max\limits_{1\leq k\leq d}\left(\mathsf{Re}\left(\lambda_{k}\right)+\sum_{j=1}^{L}
(\mathsf{Re}(\mu_{k,j}))^2\right)<0.
\end{equation}
For any positive increasing sequence $(s_n)_{n\in \mathbb{N}}$ such that $s_n\to \infty$ as $n\to \infty$, for any positive non-increasing sequence $\overline \e:= (\e_n)_{n\in \mathbb{N}}$, $x\in \mathbb{R}^d$, $x\neq 0_d$, $p>1$ it follows for any $n\in \mathbb{N}$
\begin{equation}
\begin{split}
\mathbb{P}\left(\sup_{t\in [s_n, s_{n+1}]} \|X_t(x)\| > \e_n \|x\|\right)&\leq \frac{1}{\e^2_n \|x\|^2}
\left(\frac{p}{p-1}\right)^p
\exp\left(2(p^2-p)\delta_{n}\max\limits_{1\leq k\leq d}\sum_{j=1}^L (\mathsf{Re}(\mu_{k,j}))^2\right)\\
&\qquad\times
\sum_{k=1}^{d} |\langle x, v_k\rangle|^2  \exp\Big((2 \mathsf{Re}(\widetilde{\lambda}_k)  +  \sum_{j=1}^L |\mu_{k,j}|^2)s_n\Big),
\end{split}
\end{equation}
where $\delta_n:=s_{n+1}-s_n>0$.
In particular, for 
\begin{equation}\label{e:uniformLambdauno}
\Lambda(x) :=\max\limits_{\substack{1\leq k\leq d\\ \langle x, v_k\rangle \neq 0}}\left(\mathsf{Re}\left(\lambda_{k}\right)+\sum_{j=1}^{L}
(\mathsf{Re}(\mu_{k,j}))^2\right)<0
\end{equation}
 it follows that
\begin{equation}\label{e:expprobdecay}
\begin{split}
&\mathbb{P}\left(\sup_{t\in [s_n, s_{n+1}]} \|X_t(x)\| > \e_n \|x\|\right)\\
&\quad\leq \frac{\exp\Big(2\Lambda(x) s_n\Big)}{\e^2_n}
\left(\frac{p}{p-1}\right)^p
\exp\left(2(p^2-p)\delta_{n}\max\limits_{1\leq k\leq d}\sum_{j=1}^L (\mathsf{Re}(\mu_{k,j}))^2\right)
  ,\qquad n\in \mathbb{N}.
\end{split}
\end{equation}
\end{proposition}
The proof is given in Section~\ref{a.s.}.

Note that for generic values $x\in \RR^d$ we have $\Lambda(x) = \Lambda$, that is, the right-hand side does not depend on $x$.

With the help of the following quantified version of the first Borel--Cantelli Lemma, the preceding rates can be translated into the following trade-off between a.s. rates and the tail decays of the (random) modulus of convergence.

\begin{lemma}[Quantified version of the first Borel–Cantelli Lemma, see~Lemma~2~in~\cite{EHS}]\label{lem:BC}\hfill

\noindent 
Consider a sequence of random vectors $(\mathcal{Z}_n)_{n\in \mathbb{N}}$ and $\mathcal{Z}$ with values in $\RR^d$ on a given probability space $(\Omega, \mathcal{F}, \mathbb{P})$. Denote for any $\e>0$ the anti-concentration function  
\[
p(n;\e) := \mathbb{P}(\|\mathcal{Z}_n - \mathcal{Z}\|>\e).
\]
Assume that for any fixed $\e>0$ the convergence in probability, that is, 
$p(n;\e)\to 0$, as $n\to \infty$. 
Then for any positive non-increasing sequence $\overline \e := (\e_n)_{n\in \mathbb{N}}$ and any positive non-decreasing sequence $\overline w:= (w_n)_{n\in \mathbb{N}}$ such that 
\begin{equation}
K(\overline w;\overline \e) := \sum_{n=0}^\infty w_n \sum_{m=n}^\infty p(m; \e_m) <\infty
\end{equation}
the $\overline \e$-(random) modulus of convergence 
\[
\mathcal{M}_{\overline \e}:= \max\{n\in \mathbb{N}~|~\|\mathcal{Z}_n-\mathcal{Z}\|>\e_n\}
\]
satisfies the following statements:
\begin{enumerate}
\item $\|Z_n -Z\|  \lqq \e_n$ for all $n > \mathcal{M}_{\overline \e}$ a.s. 
\item For $F(n):= \sum\limits_{j=0}^{n-1} w_j$, $n\in \mathbb{N}$ with $F(0) = 0$ it follows that
\begin{equation}
\mathbb{E}[F(\mathcal{M}_{\overline \e})] \leq K(\overline w; \overline \e).
\end{equation}
In particular, $\mathbb{P}(\mathcal{M}_{\overline \e}\gqq r)\leq (F(r))^{-1} K(\overline w; \overline \e)$ for all $r\in \mathbb{N}$. 
\end{enumerate}
\end{lemma}

\begin{remark}[The trade-off relation between modulus of convergence and almost sure upper bounds]\hfill

\noindent
The preceding lemma is a quantification of the well-known fact that the summability of the probability of the error event (anti-concentration probabilities $\{\|Z_n-Z\|>\e\}$ for any fixed $\e>0$) implies a.s. convergence. 

It is worth noting the trade-off relation between the sequence $\overline \e = (\e_n)_{n\in \NN}$ and $\overline w=(w_n)_{n\in \NN}$. The faster $\overline \e = (\e_n)_{n\in \mathbb{N}}$ tends to $0$ as a function of $n$ the slower the sequence $p(n; \e_n)$ tends to $0$ (as a function of $n$) and hence the slower $w_n$ (and its discrete antiderivative $F(n)$, respectively) can tend to $\infty$, which finally results in moments of lower order for $\mathcal{M}_{\overline \e}$. 
This trade-off is optimal up to one polynomial order of integrability of $\mathcal{M}_{\overline \e}$. 
For an in-depth discussion, the proofs and several classes of examples we refer to~\cite{EHS}.
\end{remark}

We implement the trade-off for the a.s. convergence of $X_t(x)$ to $0_d$ as $t\to \infty$. 

\begin{corollary}[The a.s. trade-off for multivariate geometric Brownian motion]\label{cor:BCsupremum} 
\hfill

\noindent
Let the hypotheses of Proposition~\ref{thm:BCsupremum} be satisfied for $x\in \RR^d$, $x\neq 0_d$.
Then 
\begin{itemize}
\item[(a)]
for any positive non-increasing sequence $\overline \e:= (\e_n)_{n\in \mathbb{N}}$, 
\item[(b)] any 
positive non-decreasing sequence $\overline w:= (w_n)_{n\in \mathbb{N}}$ 
\item[(c)] any positive strictly increasing diverging sequence $\overline s := (s_n)_{n\in \NN}$ with $s_0 = 0$,
\item[(d)]
and
$\overline \delta_n := (\delta_n)_{n\in \NN}$ given by
$\delta_n := s_{n} -s_{n-1} $, $n\in \mathbb{N}$,
\end{itemize}
satisfying 
for any $p>1$ fixed
that the following constant $K:=K_{\overline \e, \overline s, \overline w, \eta,p}(x)$ is finite, where
\begin{equation}\label{eq:defKf}
\begin{split}
K:=\Big(\frac{p}{p-1}\Big)^p\sum_{n=0}^{\infty} w_n 
\sum_{m=n}^\infty
\frac{\exp\Big(2 \Lambda(x) s_m\Big)}{\e^2_m} 
\exp\left(2\delta_{m}(p^2-p)\max\limits_{1\leq k\leq d}\sum_{j=1}^L (\mathsf{Re}(\mu_{k,j}))^2\right), 
\end{split}
\end{equation}
we have for
\begin{equation}
\mathcal{M}_{\overline \e, \overline s, x}(\omega) := \max\left\{n\in \NN~\Big|~\sup_{t\in [s_n, s_{n+1}]} \|X_{t}(x, \omega)\| > \e_n \|x\|\right\}
\end{equation}
the following statements: 
\begin{enumerate}
    \item  $n> \mathcal{M}_{\overline \e, \overline s, x}$ implies that
\begin{equation}
\sup_{t\in [s_n, s_{n+1}]} \|X_{t}(x)\| \lqq \e_n \|x\| \qquad \textrm{ a.s. }
\end{equation}
\item for some $p>1$ it follows that
\begin{equation}
\mathbb{E}\big[F(\mathcal{M}_{\overline \e, \overline s, x})\big] \leq 
 K,
\end{equation}
in particular, 
\[
\mathbb{P}(\mathcal{M}_{\overline \e, \overline s, x}\gqq r) \lqq \min\left\{\frac{K}{F(r)},1\right\},\quad r\in \mathbb{N}.
\]
\end{enumerate}
\end{corollary}
The proof is a direct consequence of Lemma~\ref{lem:BC} with the help of Proposition~\ref{thm:BCsupremum}.

In the sequel, we explain two extreme cases. We first see that the propagation along a natural time scale 
on intervals of order $(n, n+1]$ yields asymptotically a.s. exponential rates with exponential integrability of the random modulus of convergence. 

\begin{example}[Almost sure exponential rates of convergence]\hfill\\
Consider the hypotheses of Corollary~\ref{cor:BCsupremum}
with $s_n = n$ for all $n\in \mathbb{N}$.
Then for any $p>1$
and any positive non-increasing sequence $\overline \e:= (\e_n)_{n\in \mathbb{N}}$
it follows that
\begin{equation}
\begin{split}
\mathbb{P}\left(\sup_{t\in (n, n+1]} \|X_t(x)\| > \e_n \|x\|\right)\leq  \frac{\pi_n}{\e^2_n},
\end{split}
\end{equation}
where 
\[
\pi_n:=C_{p,L,d} ~\exp\Big(-2 |\Lambda(x)| n\Big),\quad n\in \mathbb{N},
\]
and
\[
C_{p,L,d}:=\Big(\frac{p}{p-1}\Big)^p ~\exp\Big(2p(p-1)  \max\limits_{1\leq k\leq d}\sum_{j=1}^L (\mathsf{Re}(\mu_{k,j}))^2\Big).
\]
Assume that~\eqref{eq:defKf}
is finite for a suitable positive non-decreasing sequence $\overline w:= (w_n)_{n\in \mathbb{N}}$.
Then for the choice  $\e_n:= e^{-r |\Lambda(x)| n}$, $n\in \mathbb{N}$ with $r\in (0,1)$ 
it follows that
\begin{itemize}
\item[(i)] for $n> \mathcal{M}_{\overline \e, \overline s, x}$ we have 
\begin{equation}
\sup_{t\in (n,  
n+1 ]} \|X_t(x)\| \leq e^{-r |\Lambda(x)| n} \|x\| \qquad \textrm{ a.s.},
\end{equation}
\item[(ii)]
for $w_n:=e^{2q(1-r) |\Lambda(x)|n}$, $n\in \mathbb{N}$,
with $q\in (0, 1)$ we have 
\[
\EE\Big[e^{2q(1-r) |\Lambda(x)|\mathcal{M}_{\overline \e, \overline s, x}}\Big] \lqq 1+ \sum_{n=1}^\infty w_n\sum_{m=n}^\infty 
\min\left\{\frac{\pi_m}{\e^2_m},1\right\}<\infty.
\]
\end{itemize}
The preceding upper bounds can be calculated explicitly. Note $r \in (0,1)$ can be chosen in a way that 
$q$ can be chosen arbitrarily close to $1$. Here, you can see the trade-off between the asymptotic upper bound for the a.s. rate and the modulus of convergence. The closer $r$ is chosen to be to $1$, the lower is the exponent $q$, of the modulus and vice-versa. 
\end{example}

\begin{example}[The trade-off for almost sure upper pre-cutoff convergence]\label{ex:precutoff}\hfill\\
Consider the hypotheses of Corollary~\ref{cor:BCsupremum} for some $x\in \RR^d$, $x\neq 0_d$, 
with running time
\[s_n =(1+c) \frac{1}{|\Lambda(x)|} \ln(n+1)\quad \textrm{ for all } \quad n\in \mathbb{N},
\]
where $c>1$ is fixed.
This implies
$\delta_n\leq C$ for all $n\in \mathbb{N}$ with  $C:=(1+c)|\Lambda(x)|^{-1}>0$. 
Hence, the following a.s. trade-off for the upper pre-cutoff convergence is valid. 
Whenever the corresponding expression in~\eqref{eq:defKf} is finite, we have a.s. for all $n> \mathcal{M}_{\overline \e, \overline s, x}$
\begin{equation}
\sup_{t\in (s_n, s_{n+1}]} \|X_t(x)\| \lqq \e_n \|x\|,\quad \e_n>0,
\end{equation}
and for any positive, non-decreasing sequence $\overline w:= (w_n)_{n\in \NN}$ such that 
\begin{equation}
\begin{split}
K:=C_{p,L,d}
\sum_{n=0}^{\infty} w_n 
\sum_{m=n}^\infty
\frac{1}{\e^2_m}
\frac{1}{(m+1)^{2+2c}}<\infty 
\end{split}
\end{equation}
where
\[
C_{p,L,d}:=\Big(\frac{p}{p-1}\Big)^p
\exp\left(2C(p^2-p)\max\limits_{1\leq k\leq d}\sum_{j=1}^L (\mathsf{Re}(\mu_{k,j}))^2\right).
\]
For the choice $\e_n:= \frac{1}{n+1}$ and $w_n=n^{r}$ for some $r\geq 0$ and $n\in \mathbb{N}$ we have
\begin{equation}
\begin{split}
\sum_{n=0}^{\infty} w_n 
\sum_{m=n}^\infty
\frac{1}{(m+1)^{2c}}<\infty\quad \textrm{ when }\quad c>1+r/2.
\end{split}
\end{equation}
\end{example}
\bigskip 
\section{\textbf{Proof of: Cutoff convergence for the autocorrelation function
}}\label{sec:thprauto} 

In this section, we prove Theorem~\ref{thm:auto}, Theorem~\ref{thm:autoprof} and Corollary~\ref{cor:norm}.
\begin{proof}[Proof of Theorem~\ref{thm:auto}]
Recall the definition of
$t_\e$ in~\eqref{eq:mixtime}.
Straightforward computations gives
\begin{equation}
\lim\limits_{\e\to 0}\frac{(t_\e+r\cdot w_\e)^{{\ell}-1}e^{-{q} (t_\e+r\cdot w_\e)}}{\e}={q}^{1-{\ell}}e^{-r{q}w}\quad \textrm{ for all }\quad r\in \mathbb{R}.
\end{equation}
By~\eqref{eq:normrelation} we have
\begin{equation}\label{eq:HSt}
\begin{split}
\frac{1}{\e}\left\|R_{t_\e+r\cdot w_\e}(x)\right\|_{\mathrm{HS}}
=\frac{(t_\e+r\cdot w_\e)^{{\ell}-1}} {e^{{q} (t_\e+r\cdot w_\e)}}\frac{1}{\e}\left\| \frac{e^{{q} (t_\e+r\cdot w_\e)}}{(t_\e+r\cdot w_\e)^{{\ell}-1}} e^{\Delta (t_\e+r\cdot w_\e)}\mathrm{vec}(xx^*)\right\|_{\mathbb{R}^{d^2}}.
\end{split}
\end{equation}
Then~\eqref{eq:belowaboveti},~\eqref{eq:limiteinf}~and~\eqref{eq:limitesup} give for all $r\in \mathbb{R}$
\begin{equation}
\begin{split}
\liminf\limits_{\e\to 0}\frac{1}{\e}\left\|R_{t_\e+r\cdot w_\e}(x)\right\|_{\mathrm{HS}}
={q}^{1-{\ell}}e^{-r{q}w}\,\liminf\limits_{\e\to 0}
\left\|\sum_{k=1}^{{m}} e^{\ii  (t_\e+r\cdot w_\e){\theta}_k} {v}_k\right\|\geq 
{K}_0 {q}^{1-{\ell}}e^{-r{q}w}
\end{split}
\end{equation}
and
\begin{equation}
\begin{split}
\limsup\limits_{\e\to 0}\frac{1}{\e}\left\|R_{t_\e+r\cdot w_\e}(x)\right\|_{\mathrm{HS}}
={q}^{1-{\ell}}e^{-r{q}w}\,\limsup\limits_{\e\to 0}
\left\|\sum_{k=1}^{{m}} e^{\ii  (t_\e+r\cdot w_\e){\theta}_k} {v}_k\right\|\leq 
{K}_1 {q}^{1-{\ell}}e^{-r{q}w},
\end{split}
\end{equation}
which implies~\eqref{ec:upper}~and~\eqref{ec:lower}, respectively.
\end{proof}

\begin{proof}[Proof of Theorem~\ref{thm:autoprof}]
Due to the Bolzano--Weierstrass Theorem and~\eqref{ec:uno1} we have that $\omega(x)$ is not the empty set and
$O$ does not belong to $\omega(x)$.
The proof of the statement follows straightforwardly from~\eqref{ec:uno1} and~\eqref{eq:HSt}.
\end{proof}

\begin{proof}[Proof of Corollary~\ref{cor:norm}]
On the one hand, 
the Cauchy--Schwarz inequality yields
\begin{equation}\label{eq:uno}
\begin{split}
\|R_t(x)\|_{\mathrm{HS}}:&=\left(\sum_{j,k=1}^d (R_t(x))^2_{j,k}\right)^{1/2}
=\left(\sum_{j,k=1}^d \left(
\mathbb{E}[(X_t(x))_j(X_t(x))_k]\right)^2\right)^{1/2}\\
&
\leq 
\left(\sum_{j,k=1}^d 
\mathbb{E}[(X_t(x))^2_j]
\mathbb{E}[(X_t(x))^2_k]\right)^{1/2}=
\EE\left[\|X_{t}(x)\|^2_{\mathbb{R}^d}\right].
\end{split}
\end{equation}
On the other hand,
\begin{equation}\label{eq:dos}
\begin{split}
\|R_t(x)\|_{\mathrm{HS}}&\geq
\left(\sum_{j=1}^d \left(
\mathbb{E}[(X_t(x))^2_j]\right)^2\right)^{1/2}
\geq \frac{1}{\sqrt{d}}
\sum_{j=1}^d 
\mathbb{E}[(X_t(x))^2_j]=\frac{1}{\sqrt{d}}\EE\left[\|X_{t}(x)\|^2_{\mathbb{R}^d}\right].
\end{split}
\end{equation}
By~\eqref{eq:uno}~and~\eqref{eq:dos} we obtain
\begin{equation}\label{eq:bothine}
\frac{1}{\sqrt{d}}\EE\left[\|X_{t}(x)\|^2_{\mathbb{R}^d}\right]
\leq \|R_t(x)\|_{\mathrm{HS}}\leq  
\EE\left[\|X_{t}(x)\|^2_{\mathbb{R}^d}\right].
\end{equation}
Therefore,~\eqref{eq:bothine}
with the help of~\eqref{ec:lower} in Theorem~\ref{thm:auto} gives
\begin{equation}
\infty=
\lim_{r \to -\infty} \liminf_{\e\to 0}
\frac{1}{\e}\|R_{t_\e+r\cdot w_\e}(x)\|_{\mathrm{HS}}\leq 
\lim_{r \to -\infty} \liminf_{\e\to 0}
\frac{1}{(\sqrt{\e})^2}
\EE\left[\|X_{t_\e+r\cdot w_\e}(x)\|^2_{\mathbb{R}^d}\right]\leq \infty,
\end{equation}
that is,
\begin{equation}
\lim_{r \to -\infty} \liminf_{\e\to 0}
\frac{1}{(\sqrt{\e})^2}
\EE\left[\|X_{t_\e+r\cdot w_\e}(x)\|^2_{\mathbb{R}^d}\right]=\infty.
\end{equation}
Also, ~\eqref{eq:bothine}
with the help of~\eqref{ec:upper} in Theorem~\ref{thm:auto} yields
\begin{equation}
0\leq \lim_{r \to \infty} \limsup_{\e\to 0}
\frac{1}{(\sqrt{\e})^2}
\EE\left[\|X_{t_\e+r\cdot w_\e}(x)\|^2_{\mathbb{R}^d}\right]\leq \sqrt{d}\cdot \lim_{r \to \infty} \limsup_{\e\to 0}
\frac{1}{\e}\|R_{t_\e+r\cdot w_\e}(x)\|_{\mathrm{HS}}=0
\end{equation}
that is,
\begin{equation}
\lim_{r \to \infty} \limsup_{\e\to 0}
\frac{1}{(\sqrt{\e})^2}
\EE\left[\|X_{t_\e+r\cdot w_\e}(x)\|^2_{\mathbb{R}^d}\right]=0.
\end{equation}
This completes the proof.
\end{proof}
 
\section{\textbf{Proof of: Cutoff convergence for simultaneously diagonalizable coefficients}}\label{sec:thpr} 

In this section, we show Lemma~\ref{lemma:diagonal}, Lemma~\ref{lem:keep}, Theorem~\ref{thm:solodiagonalizable}, Corollary~\ref{thm:profdiagonalizable}, Corollary~\ref{thm:diagonalizable} and Corollary~\ref{cor:mixingI}.

\begin{proof}[Proof of Lemma~\ref{lemma:diagonal}]
By~\eqref{eq:simdia} we have that all elements of $\{A,B_1,\ldots,B_L\}$ are pairwise commuting.
Then the solution of~\eqref{eq:mstra} is given by
\begin{equation}
X_t(x) 
= \exp\Big(A t + \sum_{j=1}^L B_j W_t^j\Big) x
\end{equation}
for all $t\geq 0$ and $x\in \mathbb{R}^d$, see  for instance, Section~3.4 of~\cite{Mao} or Section~4.8 of~\cite{Kloeden}. 
Since the set $\{A, B_1, \dots, B_L\}$ is simultaneously diagonalizable (see Definition~\ref{def:sdud}),
we have
\[
X_t(x)= P^{-1} \exp\Big(D_A t + \sum_{j=1}^LD_{B_j} W_t^j \Big) Px
\]
for all $t\geq 0$ and $x\in \mathbb{R}^d$.
Recalling that $M^* = {\overline{M}}^{T}$ for $M\in \mathbb{C}^{d\times d}$, we then obtain
\begin{equation}
\begin{split}
\|X_t(x)\|^2 
&= x^* P^* \exp\Big(D_A^* t  + \sum_{j=1}^LD_{B_j}^* W_t^j \Big) (P^{-1})^* P^{-1} \exp\Big(D_A t  + \sum_{j=1}^LD_{B_j} W_t^j \Big) P x\\
&= y^* \exp\Big(D_A^* t  + \sum_{j=1}^LD_{B_j}^* W_t^j \Big) \Gamma \exp\Big(D_A t  + \sum_{j=1}^LD_{B_j} W_t^j \Big) y,
\end{split}
\end{equation}
where 
$\Gamma:= (P^{-1})^* P^{-1}$ and 
$y:=Px$.  
It is easy to see that $\Gamma$ is a Hermitian matrix and it has non-negative eigenvalues and real-valued eigenvectors.
Since the expectation is a linear operator, we have 
\begin{equation}\label{eq:E2}
\mathbb{E}[\|X_t(x)\|^2] 
=y^* \mathbb{E}\Big[H(t)\Big] y,
\end{equation}
where 
\[
H(t):=\exp\Big(D_A^* t  + \sum_{j=1}^LD_{B_j}^* W_t^j \Big) \Gamma \exp\Big(D_A t  + \sum_{j=1}^LD_{B_j} W_t^j \Big),\quad t\geq 0.
\]
On the one hand we obtain that the random matrix $H(t)$ under the expectation can be evaluated component-wise due to the diagonal structure of the exponentials, that is, for each $i,k\in \{1,\ldots,d\}$ we have
\begin{equation}
\begin{split}
(H(t))_{i,k} &= 
\left(\exp\Big(D_A^* t  + \sum_{j=1}^LD_{B_j}^* W_t^j \Big)\right)_{i,i} \Gamma_{i,k} \left(\exp\Big(D_A t  + \sum_{j=1}^LD_{B_j} W_t^j \Big)\right)_{k,k}\\
&=\exp\Big(\overline{\alpha}_i t  + \sum_{j=1}^L \overline{\beta}_{i,j} W_t^j \Big) \Gamma_{i,k} \exp\Big(\alpha_k t  + \sum_{j=1}^L \beta_{k,j} W_t^j \Big)\\
&=\exp\Big((\overline{\alpha}_i+\alpha_k) t  + \sum_{j=1}^L (\overline{\beta}_{i,j}+\beta_{k,j}) W_t^j \Big) \Gamma_{i,k}
\end{split}
\end{equation}
and whose expectation is given by
\begin{equation}
\begin{split}
(\mathbb{E}[H(t)])_{i,k}=\mathbb{E}[(H(t))_{i,k}]
&= \Gamma_{i, k} \exp\big((\overline \alpha_i + \alpha_k) t \big) \mathbb{E}\Big[\exp\big(\sum_{j=1}^L (\overline \beta_{i, j}+\beta_{k, j}) W_t^j\big)\Big]\\
&= \Gamma_{i, k} \exp\big((\overline \alpha_i + \alpha_k) t \big) \exp\left(\frac{1}{2} \sum_{j=1}^L (\overline \beta_{i, j}+\beta_{k, j})^2 t\right).
\end{split}
\end{equation}
On the other hand, we note that for each $i,k\in \{1,\ldots,d\}$ 
\begin{equation}
(\mathbb{E}[H(t)])_{i,k} 
= \Gamma_{i, k}\exp(t R_{i,k})\exp(\mathsf{i} t I_{i,k})
\end{equation}
with
\begin{equation}
\begin{split}
R_{i,k}:&=\textsf{Re}\left(\overline{\alpha}_i + \alpha_k  +\frac{1}{2}\sum_{j=1}^L (\overline{\beta}_{i,j} + \beta_{k,j})^2
\right),\\
I_{i,k}:&=\textsf{Im}\left(\overline{\alpha}_i + \alpha_k  +\frac{1}{2}\sum_{j=1}^L (\overline{\beta}_{i,j} + \beta_{k,j})^2
\right).
\end{split}
\end{equation}
By~\eqref{eq:E2} we have the desired representation, that is,
\begin{equation}
\begin{split}
\mathbb{E}[\|X_t(x)\|^2] 
=\sum_{i,k=1}^{d}\overline{y}_i \Gamma_{i, k}y_k\exp(t R_{i,k})\exp(\mathsf{i} tI_{i,k}),\quad x\in \mathbb{R}^d,\,t\geq 0.
\end{split}
\end{equation}
\end{proof}

\begin{proof}[Proof of Lemma~\ref{lem:keep}]
By Lemma~\ref{lemma:diagonal} we have
\begin{equation}
\begin{split}
\mathbb{E}[\|X_t(x)\|^2] 
=
\sum_{i,k=1}^{d}
\overline{y}_i \Gamma_{i, k}y_k\exp(-t (-R_{i,k}))\exp(\mathsf{i} tI_{i,k})\mathbbm{1}(\overline{y}_i \Gamma_{i, k}y_k\neq 0).
\end{split}
\end{equation}
Since 
$\|x\|^2=\sum_{i,k=1}^{d}\overline{y}_i \Gamma_{i, k}y_k>0$, we have existence of $(i_0,k_0)\in \{1,\ldots,d\}^2$ (and hence $(k_0,i_0)\in \{1,\ldots,d\}^2$) such that $\overline{y}_{i_0} \Gamma_{i_0, k_0}y_{k_0}\neq 0$ (and also $\overline{y}_{k_0} \Gamma_{i_0, k_0}y_{i_0}\neq 0$).
Let 
\begin{equation}
q:=q(x)=\min\{-R_{i,k}:i,k\in \{1,\ldots,d\}\quad \textrm{ and }\quad\overline{y}_i \Gamma_{i, k}y_k\neq 0\}
\end{equation}
and define
\[
I:=I(x)=\{(i,k)\in \{1,\ldots,d\}^2: -R_{i,k}=q\quad \textrm{ and }\quad \overline{y}_i \Gamma_{i, k}y_k\neq 0\}.
\]
Clearly $I$ is not an empty set. Moreover, $I$ is a symmetric set, that is, $(i,k)\in I$ iff $(k,i)\in I$. 
Since the cardinality of $I$ is finite, 
a diagonal Cantor argument and the Bolzano–Weierstrass Theorem imply
\begin{equation}
\begin{split}
\liminf\limits_{t\to \infty}e^{qt}\mathbb{E}[\|X_t(x)\|^2] 
&=2\sum_{(i,k)\in I;\,i\leq k}\mathrm{Re}\left(\overline{y}_i \Gamma_{i, k}y_k\exp(\mathsf{i} l_{i,k})\right)\quad \textrm{ and }\\
\limsup\limits_{t\to \infty}e^{qt}\mathbb{E}[\|X_t(x)\|^2] 
&=2\sum_{(i,k)\in I;\,i\leq k}\mathrm{Re}\left(\overline{y}_i \Gamma_{i, k}y_k\exp(\mathsf{i} u_{i,k})\right)
\end{split}
\end{equation}
for some $l_{i,k},u_{i,k}\in \mathbb{R}$.
Moreover, $\liminf\limits_{t\to \infty}e^{qt}\mathbb{E}[\|X_t(x)\|^2]$ exists if and only if~\eqref{eq:profcond} is valid. This completes the proof.
\end{proof}

\begin{proof}[Proof of Theorem~\ref{thm:solodiagonalizable}]
By~\eqref{ec:limsupinff}~in~Lemma~\ref{lem:keep} we have
\begin{equation}\label{eq:liminfec1}
\begin{split}
\liminf\limits_{\e\to 0}
\frac{\mathbb{E}[\|X_{t_\e+r\cdot w_\e}(x)\|^2]}{\e}&=
\liminf\limits_{\e\to 0}\frac{e^{-q(t_\e+r\cdot w_\e)}}{\e}e^{q(t_\e+r\cdot w_\e)}\mathbb{E}[\|X_{t_\e+r\cdot w_\e}(x)\|^2]\\
&=\liminf\limits_{\e\to 0}e^{-qr\cdot w_\e}e^{q(t_\e+r\cdot w_\e)}\mathbb{E}[\|X_{t_\e+r\cdot w_\e}(x)\|^2]\\
&=2e^{-qr\cdot w}\sum_{(i,k)\in I;\,i\leq k}\mathrm{Re}\left(\overline{y}_i \Gamma_{i, k}y_k\exp(\mathsf{i} l_{i,k})\right)>0
\end{split}
\end{equation}
for some $l_{i,k}\in \mathbb{R}$.
Analogously,~\eqref{ec:limsupinff}~in~Lemma~\ref{lem:keep} yields
\begin{equation}\label{eq:limsupec1}
\begin{split}
\limsup\limits_{\e\to 0}
\frac{\mathbb{E}[\|X_{t_\e+r\cdot w_\e}(x)\|^2]}{\e}=
2e^{-qr\cdot w}\sum_{(i,k)\in I;\,i\leq k}\mathrm{Re}\left(\overline{y}_i \Gamma_{i, k}y_k\exp(\mathsf{i} u_{i,k})\right)<\infty
\end{split}
\end{equation}
for some $u_{i,k}\in \mathbb{R}$.
By~\eqref{eq:liminfec1} and~\eqref{eq:limsupec1} we obtain~\eqref{eq:thstaresult}.
\end{proof}

\begin{proof}[Proof of Corollary~\ref{thm:profdiagonalizable}]
The proof follows directly from~\eqref{eq:profcond} in Lemma~\ref{lem:keep}.
\end{proof}

\begin{proof}[Proof of Corollary~\ref{thm:diagonalizable}]
By~\eqref{eq:limitfp}~in~Corollary~\ref{cor:simuluni} we have   
\begin{equation}
\begin{split}
\lim\limits_{\e\to 0}\frac{\mathbb{E}[\|X_{t_\e+r\cdot w_\e}(x)\|^2]}{\e}&=\lim\limits_{\e\to 0}
\frac{e^{-q(t_\e+r\cdot w_\e)}}{\e}
e^{q(t_\e+r\cdot w_\e)}\mathbb{E}[\|X_{t_\e+r\cdot w_\e}(x)\|^2]\\
&=
e^{-qr\cdot w}\sum_{i\in I} |y_i|^2>0,
\end{split}
\end{equation}
which allows us to conclude.
\end{proof}

\begin{proof}[Proof of Corollary~\ref{cor:mixingI}]
We start by noticing that
Theorem~\ref{thm:solodiagonalizable} implies
\begin{equation}\label{eq:corte}
\lim\limits_{\e\to 0}\frac{\mathbb{E}[\|X_{c\cdot t_\e}(x)\|^2]}{\e}=
\begin{cases}
\infty & \quad\mathrm{ if }\quad c\in (0,1),\\
0      & \quad \mathrm{ if }\quad   c>1.
\end{cases}
\end{equation}
Let $\delta>0$ be fixed and choose $c>1$.
Then~\eqref{eq:corte} yields the existence of $\e_0:=\e_0(\delta,c)>0$ such that for all $\e\in (0,\e_0)$ it follows
\[
\frac{\mathbb{E}[\|X_{c\cdot t_\e}(x)\|^2]}{\e}\lqq \delta.
\]
By~\eqref{eq:monotona} we infer $\tau^x_\e(\delta)\lqq c\cdot t_\e$ for all 
 $\e\in (0,\e_0)$.
Conversely, for any $M>\delta$
and
$1/c\in (0,1)$ there exists $\e_1:=\e_1(M,\delta,c)>0$ such that for all $\e\in (0,\e_1)$ it follows
\[
\frac{\mathbb{E}[\|X_{(1/c)\cdot t_\e}(x)\|^2]}{\e}>M-\delta
\]
and by~\eqref{eq:monotona} we infer $ (1/c)\cdot t_\e\lqq \tau^x_\e(M-\delta)$ for all 
 $\e\in (0,\e_1)$.
 Therefore, we have
 \[
\limsup\limits_{\e\to 0} \frac{\tau^x_\e(\delta)}{\tau^x_\e(M-\delta)}\lqq c^2
\]
for all $c>0$.
Sending $c\to 1$ we obtain the upper bounds in~\eqref{eq:mezcla}. The lower bounds follow similarly.
\end{proof}

\bigskip 

\section{\textbf{Proof of: A.s. uniform rates and a.s. upper pre-cutoff for simultaneously unitarily diagonalizable coefficients}}\label{a.s.}

In this section, we prove Proposition~\ref{thm:BCsupremum}.

\begin{proof}[Proof of Theorem~\ref{thm:BCsupremum}]
We start with the following preliminary consideration. 
Let $\alpha,\beta$ be given complex numbers. 
Write $\alpha = \alpha_1 +\ii \alpha_2$ and $\beta = \beta_1 +\ii \beta_2$ where $\alpha_1,\alpha_2,\beta_1,\beta_2$ are real numbers.
For a standard one-dimensional real Brownian motion $W=(W_t)_{t\gqq 0}$ we consider a complex geometric Brownian motion with complex coefficients $\alpha$ and $\beta$. More precisely, we consider the
unique strong solution of the following linear SDE
\begin{equation}
\ud Y_t(y) = \alpha Y_t(y) \ud t + \beta Y_t(y) \ud W_t\quad t\gqq 0,\quad Y_0 = y\in \mathbb{C}. 
\end{equation}
It is clear that the process $Y(y):=(Y_t(y))_{t\gqq 0}$ takes values on $\mathbb{C}$.
Then,  the complex conjugate process $(\overline{Y_t(y)})_{t\gqq 0}$ solves 
\begin{equation}
\ud \overline{Y_t(y)} = \overline{\alpha} \overline{Y_t(y)} \ud t + \overline{\beta} \overline{Y_t(y)} \ud W_t\quad t\gqq 0,\quad \overline{Y_0} = \overline{y}.
\end{equation}
Here we use the fact that 
 $W=(W_t)_{t\gqq 0}$ is a standard real one-dimensional Brownian motion.
By It\^o's product rule (where $([\cdot , \cdot]_t)_{t\geq 0}$ denotes the quadratic variation process, see~\cite{Pr04}, p. 66) we have 
\begin{equation}\label{eq:SDEBS}
\begin{split}
\ud |Y_t(y)|^2 &= \ud (Y_t(y) \overline{Y_t(y)}) = Y_t(y) \ud  \overline{Y_t(y)} + \overline{Y_t(y)} \ud Y_t(y) + \ud[Y(y), \overline{Y(y)}]_t\\
&= Y_t(y)(\overline{\alpha} \overline{Y_t(y)} \ud t + \overline{\beta} \overline{Y_t(y)} \ud W_t) + \overline{Y_t(y)}\left(\alpha Y_t(y) \ud t + \beta Y_t(y) \ud W_t\right) 
+  \beta \overline \beta |Y_t(y)|^2 \ud t\\
&= \overline{\alpha} |Y_t(y)|^2 \ud t + \overline{\beta} |Y_t(y)|^2 \ud W_t + \alpha |Y_t(y)|^2 + \beta |Y_t(y)|^2 \ud W_t +  |\beta|^2 |Y_t(y)|^2 \ud t\\
&= \Big(2 \mathsf{Re}(\alpha)  +  |\beta|^2\Big)|Y_t(y)|^2 \ud t + 2 \mathsf{Re}(\beta) |Y_t(y)|^2 \ud W_t.
\end{split}
\end{equation}
After relabelling $Z_t(y)=|Y_t(y)|^2$, $t\gqq 0$, 
the SDE~\eqref{eq:SDEBS} becomes a geometric Brownian motion with real coefficients and initial datum $Z_0(y)=|y|^2$, whose explicit solution for all $t\gqq 0$ is given by
\begin{equation}\label{e: explicitly}
\begin{split}
|Y_t(y)|^2 &= |y|^2  \exp\Big((2 \mathsf{Re}(\alpha)  + |\beta|^2)t\Big) 
\exp\Big(2 \mathsf{Re}(\beta) W_t - \frac{1}{2} (2\mathsf{Re}(\beta))^2 t\Big).
\end{split}
\end{equation}
For a fixed $k\in \{1,2,\ldots,d\}$
we consider the Fourier coefficient process $F^k_t(x) := \langle X_t(x), v_k\rangle $, $t\gqq 0$. 
Recall that $v_k$ is a deterministic vector in $\mathbb{C}^d$.
Then It\^o's formula with the help of~\eqref{eq:basis} yields 
\begin{equation}
\begin{split}
\ud F^{k}_t(x) &=\langle \ud X_t(x), v_k\rangle= 
\langle \widetilde A X_t(x), v_k\rangle \ud t + \sum_{j=1}^L  \langle B_j X_t(x), v_k\rangle \ud W^j_t\\
&= \langle X_t(x), \widetilde A^* v_k\rangle \ud t + \sum_{j=1}^L  \langle X_t(x), B_j^* v_k\rangle \ud W^j_t\\
&= \widetilde{\lambda}_k F^k_t(x) \ud t + \sum_{j=1}^L \mu_{k,j} F^k_t(x) \ud W^j_t,
\end{split}
\end{equation}
where in the last equality we have used~\eqref{eq:basis}~and~\eqref{eq:veryeigen}.
In the light of~\eqref{e: explicitly} for all $t\geq 0$ we obtain 
\begin{equation}
\begin{split}
|F^k_t(x)|^2 
&=  
|\langle x, v_k\rangle|^2  \exp\Big((2 \mathsf{Re}(\widetilde{\lambda}_k)  +  \sum_{j=1}^L |\mu_{k,j}|^2)t\Big)\\
&\qquad\times 
\exp\left( \sum_{j=1}^L 2 \mathsf{Re}(\mu_{k,j}) W_t^j  - \frac{1}{2} \sum_{j=1}^L(2\mathsf{Re}(\mu_{k,j}))^2 t)\right).
\end{split}
\end{equation}
Recall that $\{v_1, \dots, v_d\}$ is an orthonormal basis of $\mathbb{C}^d$. Then we have 
\[
X_t(x)=\sum_{k=1}^{d} F^k_t(x)v_k,\qquad t\geq 0,\, x\in \mathbb{R}^d.
\]
By Parseval--Plancherel isometry we have 
\begin{equation}
\begin{split}
\|X_t(x)\|^2 & =\sum_{k=1}^{d} |F^k_t(x)|^2\\
&=
\sum_{k=1}^{d} |\langle x, v_k\rangle|^2  \exp\Big((2 \mathsf{Re}(\widetilde{\lambda}_k)  +  \sum_{j=1}^L |\mu_{k,j}|^2)t\Big)M^k_t
\end{split}
\end{equation}
for all $t\geq 0,\, x\in \mathbb{R}^d$, where
\[
M^k_t:=\exp\left( \sum_{j=1}^L 2 \mathsf{Re}(\mu_{k,j}) W_t^j  - \frac{1}{2} \sum_{j=1}^L(2\mathsf{Re}(\mu_{k,j}))^2 t)\right).
\]
For each $k\in \{1,\ldots,d\}$ 
the exponential process $(M^k_t)_{t\geq 0}$ is a martingale with respect to the natural filtration of the Brownian motion $(W^k_t)_{t\geq 0}$.
We then obtain the following a priori estimate
\begin{equation}
\begin{split}
\mathbb{E}[\|X_t(x)\|^2]
&=
\sum_{k=1}^{d} |\langle x, v_k\rangle|^2  \exp\Big((2 \mathsf{Re}(\widetilde{\lambda}_k)  +  \sum_{j=1}^L |\mu_{k,j}|^2)t\Big)
\end{split}
\end{equation}
for all $t\geq 0,\, x\in \mathbb{R}^d$.

We point out that in the preceding arguments we can replace $x\in \mathbb{R}^d$  for random vectors $X\in \mathbb{R}^d$ which are independent of $\sigma(W^j_s,\,j\in \{1,\ldots,L\},\,s\geq 0)$.
The Markov property gives 
$X_{t+t'}(x)=X_{t}(X_{t'}(x))$ for all $x\in \mathbb{R}^d$, $t,t'\geq 0$.
\begin{equation}
\begin{split}
\|X_{t+t'}(x)\|^2=\|X_{t}(X_{t'}(x))\|^2
&=
\sum_{k=1}^{d} |\langle X_{t'}(x), v_k\rangle|^2  \exp\Big((2 \mathsf{Re}(\widetilde{\lambda}_k)  +  \sum_{j=1}^L |\mu_{k,j}|^2)t\Big)\\
&\qquad \times \exp\left( \sum_{j=1}^L 2 \mathsf{Re}(\mu_{k,j}) (W_{t+t'}^j-W_{t'}^j)  - \frac{1}{2} \sum_{j=1}^L(2\mathsf{Re}(\mu_{k,j}))^2 t)\right).
\end{split}
\end{equation}
By~\eqref{eq:ghy} and~\ref{e:uniformLambda} we have
\[
2 \mathsf{Re}(\widetilde{\lambda}_k)  +  \sum_{j=1}^L |\mu_{k,j}|^2<0\quad \textrm{ for all }\quad k\in \{1,\ldots,d\}.
\]
For any strictly increasing sequence $(s_n)_{n\in \NN}$ and $\lim_{n\to \infty} s_n = \infty$. 
Define $\delta_n:=s_{n+1}-s_n$, $n\in \mathbb{N}$.
Then we have
\begin{equation}
\begin{split}
\sup_{0\leq t\leq \delta_n}\|X_{t+s_{n}}(x)\|^2
&\leq 
\sum_{k=1}^{d} |\langle X_{s_{n}}(x), v_k\rangle|^2  \sup_{0\leq t\leq \delta_n}\exp\Big((2 \mathsf{Re}(\widetilde{\lambda}_k)  +  \sum_{j=1}^L |\mu_{k,j}|^2)t\Big)\\
&\qquad \times \sup_{0\leq t\leq \delta_n}\exp\left( \sum_{j=1}^L 2 \mathsf{Re}(\mu_{k,j}) (W_{t+s_n}^j-W_{s_n}^j)  - \frac{1}{2} \sum_{j=1}^L(2\mathsf{Re}(\mu_{k,j}))^2 t)\right)\\
&\leq 
\sum_{k=1}^{d} |\langle X_{s_{n}}(x), v_k\rangle|^2\\
&\qquad \times \sup_{0\leq t\leq \delta_n}\exp\left( \sum_{j=1}^L 2 \mathsf{Re}(\mu_{k,j}) (W_{t+s_n}^j-W_{s_n}^j)  - \frac{1}{2} \sum_{j=1}^L(2\mathsf{Re}(\mu_{k,j}))^2 t)\right).
\end{split}
\end{equation}
Taking expectation in both sides we obtain
\begin{equation}
\begin{split}
&\mathbb{E}\left[\sup_{0\leq t\leq \delta_n}\|X_{t+s_{n}}(x)\|^2\right]\leq 
\sum_{k=1}^{d} \mathbb{E}\left[|\langle X_{s_{n}}(x), v_k\rangle|^2\right]\\
&\qquad\quad
\times 
\mathbb{E}
\left[\sup_{0\leq t\leq \delta_n}\exp\left( \sum_{j=1}^L 2 \mathsf{Re}(\mu_{k,j}) (W_{t+s_n}^j-W_{s_n}^j)  - \frac{1}{2} \sum_{j=1}^L(2\mathsf{Re}(\mu_{k,j}))^2 t)\right)\right].
\end{split}
\end{equation}
For $p>1$, Doob's maximal inequality gives
\begin{equation}
\begin{split} 
\mathbb{E}&\left[\sup_{0\leq t\leq \delta_n}\exp\left( \sum_{j=1}^L 2 \mathsf{Re}(\mu_{k,j}) (W_{t+s_n}^j-W_{s_n}^j)  - \frac{1}{2} \sum_{j=1}^L(2\mathsf{Re}(\mu_{k,j}))^2 t)\right)\right]\\
&\leq \left(\frac{p}{p-1}\right)^p
\mathbb{E}\left[\exp\left( \sum_{j=1}^L 2p \mathsf{Re}(\mu_{k,j}) (W_{\delta_n+s_n}^j-W_{s_n}^j)  - \frac{1}{2} \sum_{j=1}^L \,p(2\mathsf{Re}(\mu_{k,j}))^2 \delta_n)\right)\right]\\
&\leq \left(\frac{p}{p-1}\right)^p
\mathbb{E}\left[\exp\left( \sum_{j=1}^L 2p \mathsf{Re}(\mu_{k,j}) (W_{\delta_n+s_n}^j-W_{s_n}^j)  - \frac{1}{2} \sum_{j=1}^L \,(2p\mathsf{Re}(\mu_{k,j}))^2 \delta_n)\right)\right]\\
&\qquad \times
\exp\left( \frac{1}{2} \sum_{j=1}^L (2p\mathsf{Re}(\mu_{k,j}))^2 \delta_{n})- \frac{1}{2} \sum_{j=1}^L p(2\mathsf{Re}(\mu_{k,j}))^2 \delta_{n})\right)\\
&= \left(\frac{p}{p-1}\right)^p
\exp\left(2(p^2-p)\delta_{n}\sum_{j=1}^L (\mathsf{Re}(\mu_{k,j}))^2\right).
\end{split}
\end{equation}
Hence,
\begin{equation}
\begin{split}
&\mathbb{E}\left[\sup_{0\leq t\leq \delta_n}\|X_{t+s_{n}}(x)\|^2\right]\\
&\leq 
\sum_{k=1}^{d} \mathbb{E}\left[|\langle X_{s_{n}}(x), v_k\rangle|^2\right]\left(\frac{p}{p-1}\right)^p
\exp\left(2(p^2-p)\delta_{n}\sum_{j=1}^L (\mathsf{Re}(\mu_{k,j}))^2\right)\\
&\leq 
\left(\frac{p}{p-1}\right)^p
\max\limits_{1\leq k\leq d}
\exp\left(2(p^2-p)\delta_{n}\sum_{j=1}^L (\mathsf{Re}(\mu_{k,j}))^2\right)
\sum_{k=1}^{d} \mathbb{E}\left[|\langle X_{s_{n}}(x), v_k\rangle|^2\right]\\
&= 
\left(\frac{p}{p-1}\right)^p
\max\limits_{1\leq k\leq d}
\exp\left(2(p^2-p)\delta_{n}\sum_{j=1}^L (\mathsf{Re}(\mu_{k,j}))^2\right)
\mathbb{E}\left[
\sum_{k=1}^{d} |\langle X_{s_{n}}(x), v_k\rangle|^2\right]\\
&= 
\left(\frac{p}{p-1}\right)^p
\exp\left(2(p^2-p)\delta_{n}\max\limits_{1\leq k\leq d}\sum_{j=1}^L (\mathsf{Re}(\mu_{k,j}))^2\right)
\mathbb{E}\left[
\|X_{s_{n}}(x)\|^2\right]\\
&= 
\left(\frac{p}{p-1}\right)^p
\exp\left(2(p^2-p)\delta_{n}\max\limits_{1\leq k\leq d}\sum_{j=1}^L (\mathsf{Re}(\mu_{k,j}))^2\right)\\
&\qquad \times
\sum_{k=1}^{d} |\langle x, v_k\rangle|^2  \exp\Big((2 \mathsf{Re}(\widetilde{\lambda}_k)  +  \sum_{j=1}^L |\mu_{k,j}|^2)s_n\Big).
\end{split}
\end{equation}
Consequently, for each $\e_n>0$
the Markov inequality gives
\begin{equation}
\begin{split}
&\mathbb{P}\left(\sup_{t\in [s_n, s_{n+1}]} \|X_t(x)\| > \e_n \|x\|\right)\\
&\leq \frac{1}{\e^2_n \|x\|^2}
\left(\frac{p}{p-1}\right)^p
\exp\left(2(p^2-p)\delta_{n}\max\limits_{1\leq k\leq d}\sum_{j=1}^L (\mathsf{Re}(\mu_{k,j}))^2\right)\\
&\qquad\times
\sum_{k=1}^{d} |\langle x, v_k\rangle|^2  \exp\Big((2 \mathsf{Re}(\widetilde{\lambda}_k)  +  \sum_{j=1}^L |\mu_{k,j}|^2)s_n\Big).
\end{split}
\end{equation}
This concludes the statement.
\end{proof}

\tableofcontents

\appendix 
{
\section{ \textbf{Basic properties of the Kronecker product and their spectral properties}}\label{ap:kron}
In this section, we recall some basic properties of the Kronecker product between matrices.

\begin{lemma}[Basic properties of the Kronecker product and their spectrum]
\label{lem:kron}
\hfill

\noindent
Let $d_1,d_2,d'_1,d'_2,d''_1,d''_2\in \mathbb{N}$, $\mathbb{K}\in \{\mathbb{R},\mathbb{C}\}$ and  $\alpha_1,\alpha_2\in \mathbb{K}$.
The following statements are valid.
\begin{itemize}
\item[(i)] 

For any $M_1\in \textsf{GL}(\mathbb{K},d_1)$ and $M_2\in \textsf{GL}(\mathbb{K},d_2)$ it follows that
\[
(M_1\otimes M_2)^{-1}=M^{-1}_2\otimes M^{-1}_1.
\]
\item[(ii)]
For any $M_1 \in \mathbb{K}^{d_1\times d_1}$ 
and 
$M_2 \in \mathbb{K}^{d_2\times d_2}$
it follows that
\[
\textsf{det}(M_1\otimes M_2)=(\textsf{det}(M_1))^{d_1}(\textsf{det}(M_2))^{d_2}.
\]
\item[(iii)] For any $M_1\in \mathbb{K}^{d_1\times d'_1}$, $M_2\in \mathbb{K}^{d_2\times d'_2}$ it follows that
\[
(\alpha M_1)\otimes (\alpha_2 M_2)= \alpha_1\alpha_2 (M_1\otimes M_2).
\]
\item[(iv)] For any $M_1\in \mathbb{K}^{d_1\times d'_1}$, $M_2\in \mathbb{K}^{d_2\times d'_2}$ and $N\in \mathbb{K}^{d'_1\times d_2}$ it follows that
\[
\textsf{vec}(M_1NM_2)=(M^T_2\otimes M_1)\textsf{vec}(N).
\]
\item[(v)] 
For any $M_1\in \mathbb{K}^{d_1\times d'_1}$, $M_2\in \mathbb{K}^{d_2\times d'_2}$, $N_1\in \mathbb{K}^{d'_1\times d''_1}$
and $N_2\in \mathbb{K}^{d'_2\times d''_2}$ it follows that
\[
(M_1\otimes M_2)(N_1\otimes N_2)=(M_1N_1)\otimes (M_2N_2).
\]
\item[(vi)]
For any $M_1,M'_1\in \mathbb{K}^{d_1\times d'_1}$ and $M_2\in \mathbb{K}^{d_2\times d''_2}$ it follows that
\[
(M_1+M'_1)\otimes M_2=M_1\otimes M_2+M'_1\otimes M_2
\]
and
\[
M_2\otimes (M_1+M'_1)=M_2\otimes M_1+M_2\otimes M'_1.
\]
\end{itemize}
\end{lemma}

\bigskip

\section{\textbf{The asymptotics of the matrix exponential}}
\label{ap:exp}
The following lemma provides an asymptotic expansion of the matrix exponential for Hurwitz stable matrices.
\begin{lemma}[The asymptotics of Hurwitz stable matrix exponentials]\label{lem:jara}
\hfill

\noindent
Let $d_*\in \mathbb{N}$ be fixed. For $Q\in \RR^{d_*\times d_*}$ with $Q<0$ and any $y \in \mathbb{R}^{d_*}$, $y\neq 0_{d_*}$, there exist 
\begin{itemize}
\item[(1)] $q:=q(y,d_*)>0$,
\item[(2)] $\ell:=\ell(y,d_*) , m:=m(y,d_*)  \in  \{1,\ldots, d_*\}$,
\item[(3)] $\theta_1:=\theta_1(y,d_*),\dots,\theta_{m}:=\theta_m(y,d_*) \in \RR$,
\item[(4)] and linearly independent vectors $v_1:=v_1(y,d_*),\dots,v_{m}:=v_m(y,d_*) \in  \mathbb{C}^{d_*}$
\end{itemize}
satisfying
\begin{equation}\label{eq:vreal}
\lim_{t \to \infty} \left\|\frac{e^{q t}}{t^{\ell-1}} \exp(tQ)y - \sum_{k=1}^{m} e^{\ii  t\theta_k} v_k\right\|_{\mathbb{C}^{d_*^2}} =0.
\end{equation}
Moreover, there are positive constants $K_0 := K_0(y,d_*)$ and  $K_1 := K_1(y,d_*)$ such that 
\begin{equation}\label{eq:belowabove}
K_0\lqq \liminf_{t\rightarrow \infty}\left\|\sum_{k=1}^{m} e^{\ii  t\theta_k} v_k\right\|_{\mathbb{C}^{d_*}}
\lqq 
\limsup_{t\rightarrow \infty}\left\|\sum_{k=1}^{m} e^{\ii  t\theta_k} v_k\right\|_{\mathbb{C}^{d_*}}\lqq K_1.
\end{equation}
\end{lemma}
Since the eigenvalues of $Q\in \mathbb{R}^{d_*\times d_*}$ come in pairs of complex conjugates, the vector
$\sum_{k=1}^{m} e^{\ii  t\theta_k} v_k$ belongs to $\mathbb{R}^{d_*}$. Sometimes the preceding lemma is stated  using the Euclidean norm in $\mathbb{R}^{d_*}$.
A version of this lemma is established as 
Lemma~B.1~in~Appendix~B
in~\cite{BJ1}, p.~1195--1196, and proved there.
In~\cite{BJ1}, the result
was proven under an additional coercivity assumption for $Q$, nevertheless, a detailed examination of the proof shows that the authors make use solely of the fact that the matrix $Q<0$.

\section*{\textbf{Acknowledgments}}
\noindent
The authors are indebted with professor Jani Lukkarinen (Department of Mathematics and Statistics, University of Helsinki, Finland) for pointing out the Fuglede--Putnam--Rosenblum Theorem and a correction from the previous version manuscript. Gerardo Barrera would like to express his gratitude to the Instituto Superior T\'ecnico for all the facilities used along the realization of this work.

\noindent
\textbf{Availability of data and material.}
Data sharing not applicable to this article as no datasets were generated or analyzed during the current study.
\hfill

\noindent
\textbf{Conflict of interests.} The authors declare that they have no conflict of interest.
\hfill

\noindent
\textbf{Authors' contributions.}
All authors have contributed equally to the paper.

\noindent
\textbf{Ethical approval.} Not applicable.

\noindent
\textbf{Funding.}
Gerardo Barrera thanks the Academy of Finland, via 
the Matter and Materials Profi4 University Profiling Action, 
the Academy project No. 339228 and project No. 346306 of the Finnish Centre of Excellence 
in Randomness and STructures. 
The research of Michael A. H\"ogele was supported by the 
projects INV-2019-84-1837 and INV-2023-162-2850 of Facultad de Ciencias at Universidad de los Andes.
The research of Michael A. H\"ogele has been supported by the project ``Mean deviation frequencies and the cutoff phenomenon'' (INV-2023-162-2850) of the School of Sciences (Facultad de Ciencias) at Universidad de los Andes, Bogot\'a, Colombia.
The research of Gerardo Barrera and Michael A. H\"ogele is partially supported by European Union’s Horizon Europe research and innovation programme under the Marie Sk\l{}odowska-Curie Actions Staff Exchanges (Grant agreement No.~101183168 -- LiBERA, Call: HORIZON-MSCA-2023-SE-01).

\noindent 
\textbf{Disclaimer.} Funded by the European Union. Views and opinions expressed are however those of the author(s) only and do not necessarily reflect those of the European Union or the European Education and Culture Executive Agency (EACEA). Neither the European Union nor EACEA can be held responsible for them.

\end{document}